\documentclass[12pt]{article}

\usepackage{epsfig,epsfig}
\usepackage{epstopdf}
\usepackage{amsmath}
\usepackage{mathrsfs}
\usepackage{amsfonts}
\usepackage{bm}
\usepackage{amssymb}
\usepackage{amsthm}
\usepackage{stmaryrd}
\usepackage{booktabs}
\usepackage{color}
\usepackage{multirow}
\usepackage{graphicx}
\usepackage{subfigure}
\usepackage{float}
\usepackage{appendix}
\usepackage{tikz}
\usetikzlibrary{calc}
\usetikzlibrary{matrix}
\usepackage{mathtools}
\usepackage{algorithm,algpseudocode,float}
\usepackage{lipsum}


\usepackage{lineno}
\usepackage{srcltx,graphicx}

\newtheorem{thm}{Theorem}[section]

\newtheorem{rem}[thm]{Remark}
\newtheorem{example}[thm]{Example}

\usepackage{hyperref}

\numberwithin{equation}{section} \topmargin=-2.5cm \oddsidemargin=1cm
\begin{document}
\baselineskip=1.5pc
	
\title{
\textbf{An efficient fifth-order interpolation-based Hermite WENO scheme for hyperbolic conservation laws}
	}
	
\author{
Xiaoyang Xie\footnote{College of Engineering, Peking University, Beijing, 100871, China. Email: xiaoyangxie@stu.pku.edu.cn.},~
Zhanjing Tao\footnote{School of Mathematics, Jilin University, Changchun, 130012, China. Email: zjtao@jlu.edu.cn.},~
Chunhai Jiao\footnote{School of Mathematics, Jilin University, Changchun, 130012, China. Email: jiaoch2419@mails.jlu.edu.cn},~
and~
Min Zhang\footnote{Corresponding author. National Engineering Laboratory for Big Data Analysis and Applications, Peking University, Beijing, 100871, China.
Chongqing Research Institute of Big Data, Peking University, Chongqing, 401121, China. Email: minzhang@pku.edu.cn.}
	}
	
\date{}
\maketitle
\begin{abstract}
		In this paper, we develop a simple, efficient, and fifth-order finite difference interpolation-based Hermite WENO (HWENO-I) scheme for one- and two-dimensional hyperbolic conservation laws.
		We directly interpolate the solution and first-order derivative values and evaluate the numerical fluxes based on these interpolated values.
		We do not need the split of the flux functions when reconstructing numerical fluxes and there is no need for any additional HWENO interpolation for the modified derivative.
		The HWENO interpolation only needs to be applied one time which utilizes the same candidate stencils, Hermite interpolation polynomials, and linear/nonlinear weights for the solution and first-order derivative at the cell interface, as well as the modified derivative at the cell center.
		The HWENO-I scheme inherits the advantages of the finite difference flux-reconstruction-based HWENO-R scheme [Fan et al., Comput. Methods Appl. Mech. Engrg., 2023], including fifth-order accuracy, compact stencils, arbitrary positive linear weights, and high resolution.
		The HWENO-I scheme is simpler and more efficient than the HWENO-R scheme and the previous finite difference interpolation-based HWENO scheme [Liu and Qiu, J. Sci. Comput., 2016] which needs the split of flux functions for the stability and upwind performance for the high-order derivative terms.
		Various benchmark numerical examples are presented to demonstrate the accuracy, efficiency, high resolution, and robustness of the proposed HWENO-I scheme.
\end{abstract}
	
\vspace{5pt}
	
\noindent\textbf{Keywords:}
Interpolation-based HWENO scheme, finite difference, efficiency, hyperbolic conservation laws.
	
\newpage
	
\section{Introduction}
\label{sec:introduction}
	The high-order and high-resolution methods are important tools for simulating hyperbolic conservation laws.
	Specifically, the WENO-type schemes \cite{WENOJS1996,WENO1994,WENO-Shu2009,WENO-Shu2020,WENO-ZhuQiu2016JCP,WENO-ZhuShu2018JCP} are powerful tools due to their high-order accuracy in smooth regions and high resolution near complicated structures simultaneously escaping nonphysical oscillations adjacent to strong shocks or contact discontinuities.
	In this paper, we develop a simple and efficient fifth-order finite difference interpolation-based Hermite WENO (HWENO-I) scheme for hyperbolic conservation laws.
	
	In 1994, the WENO scheme (cell-averages version) was first introduced by Liu, Osher, and Chan \cite{WENO1994} based on essentially non-oscillatory (ENO) schemes \cite{ENO-1,ENO-2}. Later, the finite difference WENO scheme \cite{WENOJS1996} was developed by Jiang and Shu based on the flux version of the ENO scheme \cite{Shu-Osher1989} in 1996.
	Qiu and Shu developed the finite volume HWENO schemes \cite{HWENO-1,HWENO-2} following the idea of the WENO scheme \cite{WENO1994}.
	The HWENO scheme gains a more compact stencil in the spatial reconstruction than the WENO scheme with the same order of accuracy since it uses both the information of the solution and its first-order derivatives/moments.
	In \cite{HWENO-1,HWENO-2}, different stencils and reconstructed polynomials have been used for the reconstruction of solution and derivative to maintain stability.
	However, this method does not work well for strong shock problems (e.g., the double Mach and the step forward problems \cite{HWENO-2}).
	After that, several HWENO schemes have been developed. For example, Capdeville \cite{Capdeville} developed a finite volume central HWENO scheme.
	Tao et al. \cite{Tao2015JCP} constructed a finite volume central HWENO scheme on staggered meshes.
	Li et al. \cite{LiMRHW1} developed a multi-resolution HWENO scheme with unequal stencils, but the finite difference version only has a fourth-order accuracy in two dimensions.
	Liu and Qiu \cite{LiuQiu2014JSC} developed a finite difference HWENO scheme which is fifth-order in one dimension and fourth-order in two dimensions.
	Zhao and Qiu \cite{HWENO-a} developed a fifth-order finite volume HWENO scheme with artificial linear weights.
	Zhao et al. \cite{HWENO-M} proposed a fifth-order finite difference HWENO-M scheme in one and two dimensions, and modified the first-order derivative by a HWENO interpolation to improve the scheme's robustness.
	A fifth-order finite difference HWENO-L scheme \cite{HWENO-L} for one- and two- dimensional hyperbolic conservation law was developed by Zhang and Zhao, and the modified derivatives were only used in time discretization while maintaining the original derivatives in the flux reconstruction.
	Recently, Fan et al. \cite{HWENO-R} proposed a robust fifth-order finite difference HWENO (HWENO-R) scheme for compressible Euler equations. Following the idea of the HWENO-L scheme, the authors did not use the derivative value of the target cell and made the scheme simpler and more robust.
	Fan et al. \cite{HWENO-U-arxiv} proposed a fifth-order finite volume moment-based HWENO scheme with unified stencils for hyperbolic conservation laws.
	For more HWENO schemes, the interested reader can refer to \cite{HWENO-5,Tao2024JCP,7thHWENO,HWENO-3} and the references therein.

	However, the above-mentioned finite difference HWENO schemes are all based on the usual practice of directly reconstructing the flux functions for the numerical fluxes.
	{
		To achieve the linear stability of the scheme, the flux functions should be split into two parts by considering the upwinding mechanism, and then the numerical fluxes of the two parts are computed separately.}
	The upwinding of the spatial discretization provides the necessary dissipation needed in the shock-capturing scheme.
	In this paper, we present an alternative simplified interpolation-based approach to construct a fifth-order conservative finite difference HWENO scheme for one- and two-dimensional hyperbolic conservation laws, in which we directly interpolate the solution and first-order derivative values rather than the flux functions, and evaluate the numerical fluxes based on these interpolated values.

	The idea of the interpolation-based approach was first proposed in \cite{Shu-Osher1988} to construct a high-order finite difference ENO scheme.
	Then, Jiang et al. \cite{JiangShuZhang2013} constructed a finite difference WENO scheme with the Lax-Wendroff time discretization based on the idea of \cite{Shu-Osher1988}.
	Wang et al. \cite{AWENO-1} developed a fifth-order interpolation-based WENO-Z scheme for hyperbolic conservation laws to improve the accuracy of the WENO-JS scheme \cite{JiangShuZhang2013} near the critical points.
	Gao et al. \cite{AWENO-9} extended the interpolation-based WENO scheme to the seventh- and ninth-order.
	Wang et al. \cite{AWENO-MR2021} developed a fifth-order finite difference multi-resolution WENO scheme.
	Liu and Qiu \cite{LiuQiu2014JSC-2} designed a fifth-order finite difference HWENO scheme in one dimension and fourth-order in two dimensions with the interpolation-based approach.
	However, the flux splitting is also used for the approximation of high-order derivative terms in \cite{LiuQiu2014JSC-2} for the stability and the upwind performance of the scheme.
	
\begin{figure}[h]
		\centering
		\includegraphics[width=0.95\textwidth]{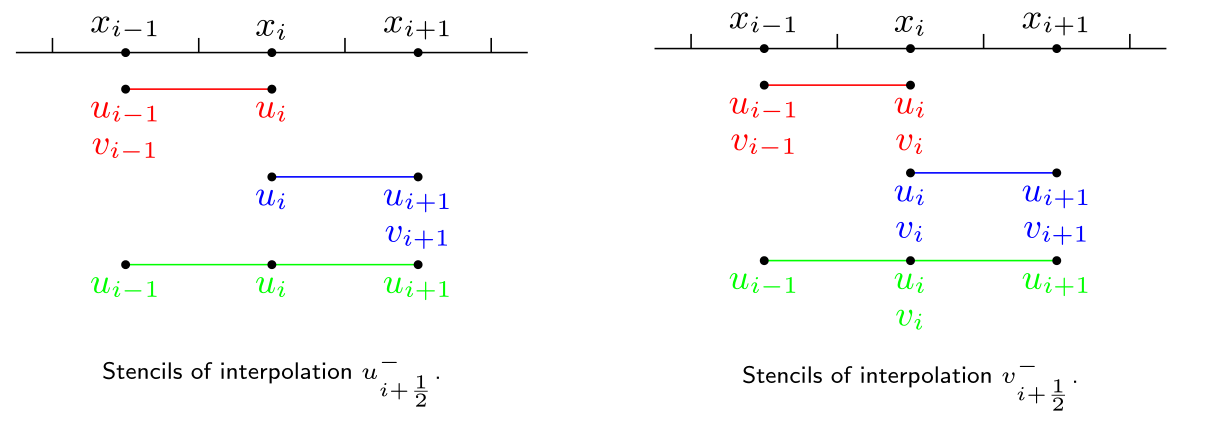}
		\caption{The stencils of interpolation for $u_{i+\frac{1}{2}}^{-}$ and $v_{i+\frac{1}{2}}^{-}$ in the previous interpolation-based {HWENO-A}  scheme \cite{LiuQiu2014JSC-2}. $v=u_x$.}
		\label{LQ-s}
	\end{figure}

\begin{figure}[h]
		\centering
		\includegraphics[width=0.5\textwidth,trim= 10 10 0 0, clip]{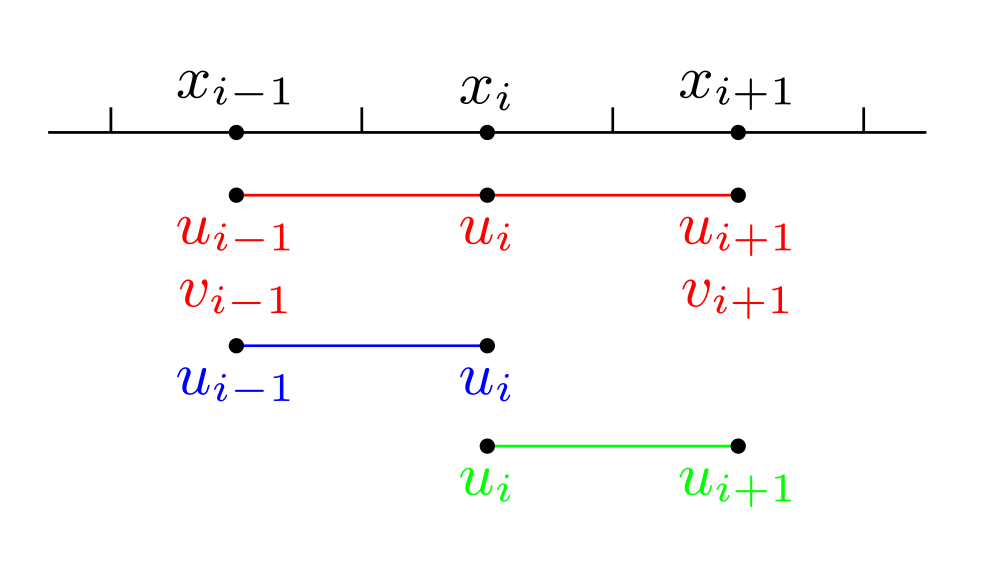}
		\caption{{
				The unified stencils of interpolation for $u_{i+\frac{1}{2}}^{-}$ and $v_{i+\frac{1}{2}}^{-}$, and modification for $v_i$, in our HWENO-I scheme. $v=u_x$.}
		}
		\label{HWENOI-s}
	\end{figure}

	Compared with the previous interpolation-based HWENO {(denoted as HWENO-A)} scheme \cite{LiuQiu2014JSC-2}, the proposed HWENO-I scheme has a much simpler implementation and better robustness.
	Specifically, our HWENO-I scheme does not need the flux splitting and directly uses the high-order central difference approximation for the high-order derivatives of the flux function.
	Second, the associated linear weights in our HWENO-I scheme are artificially set to be any positive numbers with the only requirement that their sum equals one.
	Third, our HWENO-I scheme utilizes the same candidate stencils (cf. Figs. \ref{LQ-s} and \ref{HWENOI-s}), Hermite interpolation polynomials, and linear/nonlinear weights in both interpolations of the solution and first-order derivative values which leads to better computational efficiency.
	Fourth, our HWENO-I scheme has fifth-order accuracy in both one and two dimensions, while the {HWENO-A} scheme \cite{LiuQiu2014JSC-2} only has fourth-order accuracy in two dimensions.
	Finally, our HWENO-I scheme can work well for the simulation of some strong shock problems without any PP flux limiter \cite{Hu-PP2013, XiongQiuXu-PP}, e.g., the double Mach problem (cf. Example \ref{doublemach} in \S\ref{sec:numerical}), while the {HWENO-A} scheme needs the PP limiter (cf. Remark 5 in \cite{LiuQiu2014JSC-2} ).
	In addition,  the CFL number in the {HWENO-A} scheme \cite{LiuQiu2014JSC-2} is taken to be $0.2$  while we set it as $0.6$ (cf. \S\ref{sec:numerical}) in the numerical simulations.

\begin{figure}[h]
		\centering
		\includegraphics[width=0.95\textwidth]{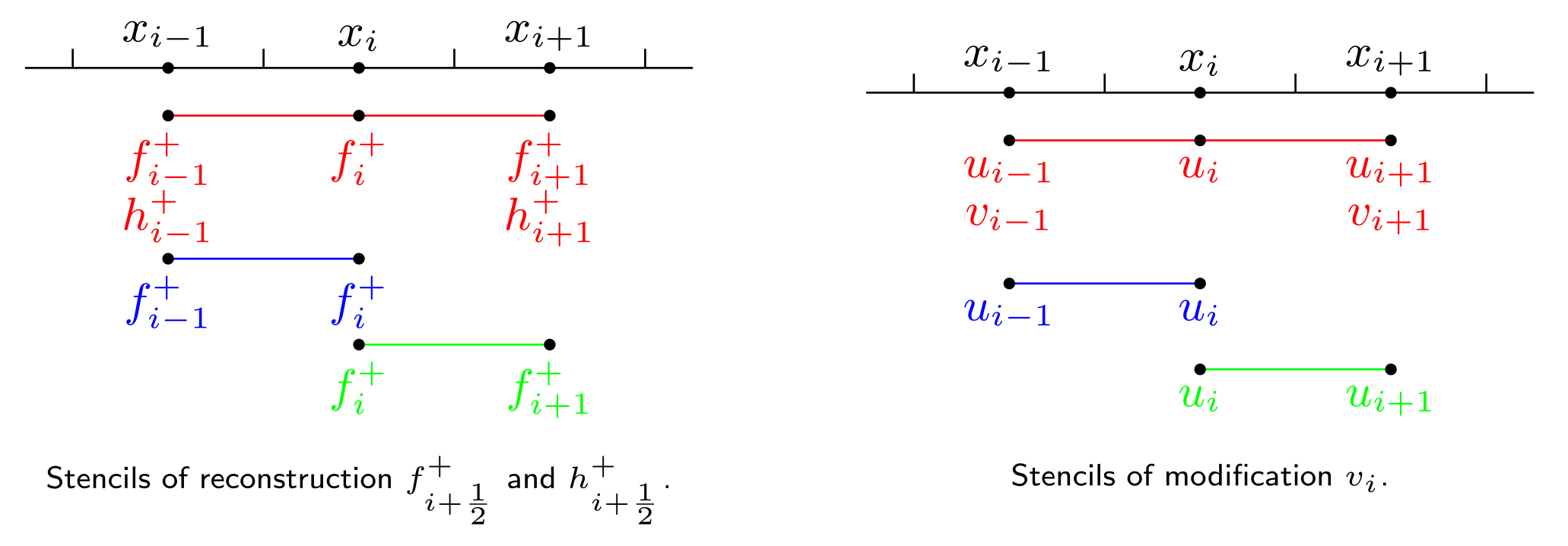}
		\caption{{
				The stencils of reconstruction for $f_{i+\frac{1}{2}}^{+}$ and $h_{i+\frac{1}{2}}^{+}$, and modification for $v_i$, in HWENO-R scheme \cite{HWENO-R}. $f=f(u)$, $h=f(u)_x$, $v=u_x$, and $df^+/du\geq0$.}}
		\label{HWENOR-s}
	\end{figure}

	Compared with the flux-reconstruction-based HWENO-R scheme \cite{HWENO-R} which has better performance in efficiency, robustness, and resolution than other finite difference HWENO schemes \cite{LiuQiu2014JSC, LiuQiu2014JSC-2, HWENO-L,HWENO-M},
	the proposed HWENO-I scheme also has some advantages.
	First, our scheme does not need the split of the flux function for the reconstruction of numerical flux.
	Second, only one HWENO interpolation is applied
	throughout the entire procedure of interpolations for the solution and first-order derivative values at the element interface, as well as the modification of the derivative value (cf. Figs.~\ref{HWENOI-s} and~\ref{HWENOR-s}).
	This leads to much simpler implementation and better computational efficiency.
	Third, the proposed HWENO-I scheme is more accurate than the HWENO-R scheme in the sense that the former can obtain a smaller error than the latter with the same mesh elements and has better efficiency and resolution for the simulation of the double Mach problem, step forward problem, etc. (cf. \S\ref{sec:numerical}).
	
	
	The organization of this paper is as follows.
	The detailed implementation of the HWENO-I scheme for one- and two-dimensional hyperbolic conservation laws are described in \S\ref{sec:HWENO} and \S\ref{sec:HWENO2d}, respectively.
	Various benchmark numerical examples (e.g., double Mach problem, step forward problem, Mach 10 shock problem with shock reflection and diffraction, Mach 2000 astrophysical jet problem, and {Rayleigh-Taylor} instability problem) are simulated in \S\ref{sec:numerical} to show the accuracy, efficiency, high resolution and robustness of the proposed HWENO-I scheme.
	Concluding remarks are given in \S\ref{sec:conclusions}.

	\section{Interpolation-based HWENO scheme in 1D}
	\label{sec:HWENO}
	
	This section presents an efficient fifth-order finite difference interpolation-based HWENO (HWENO-I) scheme for one-dimensional hyperbolic conservation laws.
	We first consider the one-dimensional scalar hyperbolic conservation laws
	\begin{equation}\label{EQ}
		\begin{cases}
			u_t + f(u)_x=0, \\
			u(x,0)=u_0(x).
		\end{cases}
	\end{equation}
	Introduce a new variable $v = u_x$, and define the new flux function $h(u,v)=f(u)_x$, i.e., $h(u,v)=f'(u)v$. {Take partial derivative with respect to $x$ on both sides of \eqref{EQ}, we can get an artificial first-order system as}
	\begin{equation}\label{EQ1}
		\begin{cases}
			u_t + f(u)_x=0, \quad u(x,0)=u_0(x),\\
			v_t + h(u,v)_x=0, \quad v(x,0)=v_0(x).
		\end{cases}
	\end{equation}
	
	Assume that the computational domain is covered by uniform meshes $I_{i} =[x_{i-\frac12},x_{i+\frac12}]$,
	~$i=1,...,N_x$,~ $x_i=\frac12( x_{i-\frac12}+x_{i+\frac12})$ with $\Delta x=x_{i+\frac12}- x_{i-\frac12}$.
	We use a high-order conservative difference scheme to discretize \eqref{EQ1} in space and get
	\begin{equation}\label{odeH}
		\begin{cases}
			\frac{d}{dt}u_i(t) =\mathcal{L}^1_{i}(u,v) := - \frac 1 {\Delta x} \big( \hat f_{i+\frac12}-\hat f_{i-\frac12}\big),\\
			\frac{d}{dt}v_i(t) =\mathcal{L}^2_{i}(u,v) := - \frac 1 {\Delta x} \big( \hat h_{i+\frac12} -\hat h_{i-\frac12}\big ),
		\end{cases}
	\end{equation}
	where $u_i(t) \approx u(x_i,t), v_i(t) \approx v(x_i,t)$, $\hat f_{i+\frac12}$ and $\hat h_{i+\frac12}$ are the numerical fluxes of the flux functions $f$ and $g$ at element interface $x_{i+\frac12}$.
	In the high-order finite difference framework, the usual practice for the reconstruction of the numerical fluxes in HWENO schemes \cite{HWENO-R,LiuQiu2014JSC,HWENO-L} is based on reconstructing the flux functions directly, and the flux function should be split by considering the upwinding mechanism.
	In this paper, we use an alternative simplified approach to construct numerical fluxes for high-order conservative finite difference schemes following the idea of \cite{LiuQiu2014JSC-2,Shu-Osher1988}.
	We define  $\hat f_{i+\frac12}$ and $\hat h_{i+\frac12}$ as follow
	\begin{equation}\label{flux-1d}
		\begin{split}
			&\hat{f}_{i+\frac12} = f(u^{-}_{i+\frac12},u^{+}_{i+\frac12})
			+\mathcal{D}_{i+\frac12}(f),\\
			&\hat{h}_{i+\frac12} = h(u^{-}_{i+\frac12},v^{-}_{i+\frac12};u^{+}_{i+\frac12},v^{+}_{i+\frac12})+\mathcal{D}_{i+\frac12}(h),
		\end{split}
	\end{equation}
	where $f(u^{-}_{i+\frac12},u^{+}_{i+\frac12})$ and $h(u^{-}_{i+\frac12},v^{-}_{i+\frac12};u^{+}_{i+\frac12},v^{+}_{i+\frac12})$ can be arbitrary monotone fluxes and
	\begin{equation}\label{flux-high}
		\begin{split}
			&\mathcal{D}_{i+\frac12}(f)\approx-\frac{1}{24}\Delta x^{2}\left(\frac{\partial^{2} f}{\partial x^{2}}\right)_{i+\frac12}
			+\frac{7}{5760}\Delta x^{4}\left(\frac{\partial^{4} f}{\partial x^{4}}\right)_{i+\frac12},\\
			&\mathcal{D}_{i+\frac12}(h)
			\approx
			-\frac{1}{24}\Delta x^{2}\left(\frac{\partial^{2} h}{\partial x^{2}}\right)_{i+\frac12}
			+\frac{7}{5760}\Delta x^{4}\left(\frac{\partial^{4} h}{\partial x^{4}}\right)_{i+\frac12}.
		\end{split}
	\end{equation}
	Here, we use the Lax-Friedrichs flux, i.e,
	\begin{equation}\label{flux-LF}
		\begin{split}
			&f(u^{-}_{i+\frac12},u^{+}_{i+\frac12})
			=\frac{1}{2}\left(f(u^{-}_{i+\frac12})+f(u^{+}_{i+\frac12})-\alpha (u^{+}_{i+\frac12}-u^{-}_{i+\frac12})\right),\quad \alpha=\max|f'(u)|,\\
			&h(u^{-}_{i+\frac12},v^{-}_{i+\frac12};u^{+}_{i+\frac12},v^{+}_{i+\frac12})=\frac{1}{2}\left(h(u^{-}_{i+\frac12},v^{-}_{i+\frac12})+h(u^{+}_{i+\frac12},v^{+}_{i+\frac12})-\alpha (v^{+}_{i+\frac12}-v^{-}_{i+\frac12})\right).
		\end{split}
	\end{equation}
	
	The scheme \eqref{odeH} has a fifth-order accuracy in space for the smooth functions $u$, $v$, $f$ and $h$ if
	\begin{equation}\label{semi-kth}
		\begin{split}
			&u^{\pm}_{i+\frac12}  = u(x_{i+\frac12},t) + \mathcal{O}(\Delta x)^5,
			\qquad
			v^{\pm}_{i+\frac12}  = v(x_{i+\frac12},t) + \mathcal{O}(\Delta x)^4,
			\\&
			\mathcal{D}_{i+\frac12}(f)=-\frac{1}{24}\Delta x^{2}\left(\frac{\partial^{2} f}{\partial x^{2}}\right)_{i+\frac12}
			+\frac{7}{5760}\Delta x^{4}\left(\frac{\partial^{4} f}{\partial x^{4}}\right)_{i+\frac12}+ \mathcal{O}(\Delta x)^5,
			\\&
			\mathcal{D}_{i+\frac12}(h)
			=
			-\frac{1}{24}\Delta x^{2}\left(\frac{\partial^{2} h}{\partial x^{2}}\right)_{i+\frac12}
			+\frac{7}{5760}\Delta x^{4}\left(\frac{\partial^{4} h}{\partial x^{4}}\right)_{i+\frac12} + \mathcal{O}(\Delta x)^{4}.
		\end{split}
	\end{equation}
	The approximations $u^{\pm}_{i+\frac12}$ and $v^{\pm}_{i+\frac12}$ are obtained by the HWENO interpolation which will be described in \S\ref{sec:HWENO-I-1d}, and the central difference approximation for
	the spatial derivatives  $\mathcal{D}_{i+\frac12}(f)$ and $\mathcal{D}_{i+\frac12}(h)$  will be described in \S\ref{sec:Dfapp}.
	
	Using the explicit third-order SSP Runge-Kutta scheme of \eqref{odeH} for time discretization, we have the fully-discrete HWENO-I scheme as, for $i=1,...,N_x,$
	\begin{equation}\label{fully-hweno}
		\begin{cases}
			\begin{cases}
				u^{(1)}_i =u^n_i+ \Delta t \mathcal{L}^1_{i}(u^{n},v^{n}),\\
				v^{(1)}_i =\tilde{v}_i ~+ \Delta t \mathcal{L}^2_{i}(u^{n},v^{n}),
			\end{cases}\\
			\begin{cases}
				u^{(2)}_i =\frac{3}{4}u^n_i+
				\frac{1}{4}\big( u^{(1)}_i+ \Delta t \mathcal{L}^1_{i}(u^{(1)},v^{(1)})\big),\\
				v^{(2)}_i =\frac{3}{4}\tilde{v}_i ~+
				\frac{1}{4}\big( \tilde{v}^{(1)}_i ~+ \Delta t  \mathcal{L}^2_{i}(u^{(1)},v^{(1)})\big),
			\end{cases}\\
			\begin{cases}
				u^{n+1}_i =\frac{1}{3}u^n_i+
				\frac{2}{3}\big( u^{(2)}_i+ \Delta t \mathcal{L}^1_{i}(u^{(2)},v^{(2)})\big),\\
				v^{n+1}_i =\frac{1}{3}\tilde{v}_i ~+
				\frac{2}{3}\big( \tilde{v}^{(2)}_i~ + \Delta t  \mathcal{L}^2_{i}(u^{(2)},v^{(2)})\big),
			\end{cases}
		\end{cases}
	\end{equation}
	where $\tilde{v}_i$, $\tilde{v}^{(1)}_i$, and $\tilde{v}^{(2)}_i$ are the fourth-order modification \eqref{limiter} of $v_i$, $v_i^{(1)}$, and  $v_i^{(2)}$, respectively.
	
	\subsection{HWENO interpolation for $u^{\pm}_{i+\frac12}$ and $v^{\pm}_{i+\frac12}$}
	\label{sec:HWENO-I-1d}
	
	Without loss of generality, we here only give the detailed reconstruction procedure for $u^{-}_{i+\frac12}$ and $ v^{-}_{i+\frac12}$, while
	the procedure for the reconstruction of $u^{+}_{i+\frac12}$ and $ v^{+}_{i+\frac12}$ is the mirror-symmetric for $x_{i+\frac12}$.
	When the reconstruction of $u^{-}_{i+\frac12}$ and $ v^{-}_{i+\frac12}$ is completed, we also accomplish the modification of the derivative value without additional HWENO interpolation.
	
	The methodology is summarized in the following.
	We choose a big stencil $T_0 = \{x_{i-1},x_i,x_{i+1}\}$
	and two small stencils $T_1 = \{x_{i-1},x_i\}$ and $T_2 = \{x_i,x_{i+1}\}$, see in Fig.~\ref{HWENOI-s}.
	There are one quintic polynomial $p_0(x)$ (on $T_0$), two linear polynomials $p_1(x)$ (on $T_1$) and $p_2(x)$ (on $T_2$), such that
	\begin{equation}\label{HWEp0}
		\begin{split}
			&p_0(x):~
			p_0(x_{i+\ell}) = u_{i+\ell},\quad \ell=-1,0,1;
			\qquad p'_0(x_{i+\ell}) = v_{i+\ell},\quad  \ell=-1,1, \\
			&p_1(x):~
			p_1(x_{i+\ell}) = u_{i+\ell},\quad \ell=-1,0, \\
			&p_2(x):~
			p_2(x_{i+\ell}) = u_{i+\ell},\quad \ell=~~~0,1.
		\end{split}
	\end{equation}
	Evaluate the values of $p_0(x),~p_1(x),~p_2(x)$ and the derivative of $p_0(x)$ at the point $x_{i+\frac12}$
	\begin{equation}\label{Pxi}
		\begin{split}			p_0(x_{i+\frac12})&=-\frac{1}{8}u_{i-1}+\frac{9}{16}u_i+\frac{9}{16}u_{i+1}
			-\Delta x \big(\frac{3}{64}v_{i-1}+\frac{9}{64}v_{i+1}\big),\\
			p_1(x_{i+\frac12})&=-\frac{1}{2} u_{i-1}+\frac{3}{2}u_i,\\
			p_2(x_{i+\frac12})&=\frac{1}{2}u_i+\frac{1}{2} u_{i+1},
		\end{split}
	\end{equation}
	and
	\begin{equation}\label{P'xi}
		\begin{split}
			p'_0(x_{i+\frac12})&=\frac{1}{\Delta x}\big(\frac{3}{16}u_{i-1}-\frac{3}{2}u_i+\frac{21}{16}u_{i+1}\big)
			+\frac{1}{16}v_{i-1}-\frac{3}{16}v_{i+1}.
		\end{split}
	\end{equation}
	Similarly, we have the derivative values of $p_0(x),~p_1(x),~p_2(x)$ at the central point $x_i$
	\begin{equation}\label{limiter-p}
		\begin{split}
			p'_0(x_i)&=
			\frac{3}{4 \Delta x}\big(u_{i+1}-u_{i-1}\big)-\frac{1}{4}\big(v_{i-1}+v_{i+1}\big),\\
			p'_1(x_i)&=\frac{1}{\Delta x}\big(u_i-u_{i-1}\big),\\
			p'_2(x_i)&=\frac{1}{\Delta x}\big(u_{i+1}-u_{i}\big),\\
		\end{split}
	\end{equation}
	which will be used in the modification of the derivative value $v_i$.

	The linear weights $\{\gamma_0,~\gamma_1,~\gamma_2\}$ can be chosen as any positive constants with $\gamma_0 +\gamma_1+\gamma_2=1$.
	To measure how smooth the functions $p_{\ell}(x),~\ell=0,1,2,$ are in the target cell $I_i$,
	we compute the smoothness indicators $\beta_{\ell}$:
	\begin{equation}
		\label{GHYZ}
		\beta_{\ell}=\sum_{\alpha=1}^k\int_{I_i}{\Delta x}^{2\alpha-1}(\frac{d ^\alpha
			p_{\ell}(x)}{d x^\alpha})^2 dx, \quad {\ell}=0,1,2,
	\end{equation}
	where $k$ is the degree of the polynomials $p_{\ell}(x)$.
	The explicit formulas are given by
	\begin{align*}
		\beta_0&=\big(a_1+\frac{1}{4}a_3\big)^2
		+\frac{13}{3}\big(a_2+\frac{63}{130}a_4\big)^2
		+\frac{781}{20}a_3^2+\frac{1421461}{2275}a_4^2,\\
		\beta_1&=(u_{i}-u_{i-1})^2, \\
		\beta_2&=(u_i -u_{i+1})^2,
	\end{align*}
	with
	\begin{align*}
		\begin{cases}
			a_1 =
			-\frac{\Delta x}{4} \big( v_{i-1}+v_{i+1}\big)
			+\frac{3}{4}\big(u_{i+1}-u_{i-1}\big),\\
			a_2 =
			\frac{\Delta x}{4} \big(v_{i-1}-v_{i+1}\big)
			+\big(u_{i-1}-2u_{i}+u_{i+1}\big),\\
			a_3 =
			\frac{\Delta x}{4} \big(v_{i-1}+v_{i+1}\big)
			+\frac{1}{4}\big(u_{i-1}-u_{i+1}\big),\\
			a_4 = \frac{\Delta x}{4} \big(v_{i+1}-v_{i-1}\big)
			-\frac{1}{2}\big(u_{i-1}-2u_{i}+u_{i+1}\big).
		\end{cases}
	\end{align*}
	Then, the nonlinear weights $\omega^{\gamma}_{\ell}$ corresponding to the linear weights $\gamma_l$ are computed as
	\begin{equation}\label{nonlinear}
		\omega^{\gamma}_{\ell}=\frac{\bar\omega_{\ell}(\gamma)}{\sum\limits_{\ell=0}^{2}\bar\omega_{\ell}(\gamma)},
		\quad \mbox{with}\quad \ \bar\omega_{\ell}(\gamma)=\gamma_{\ell}~(1+\frac{\tau}{\beta_{\ell}+\varepsilon}), \quad {\ell}=0,1,2.
	\end{equation}
	where
	$\tau=\frac{1}{4}\Big(|\beta_{0}-\beta_{1}|+|\beta_{0}-\beta_{2}|\Big)^2$, and $\varepsilon$ is a small positive number to avoid the denominator by zero. In our computation, we take $\varepsilon=10^{-10}$.
	
	Finally, the values of $u^-_{i+\frac12}$ and $v^-_{i+\frac12}$ are reconstructed by
	\begin{align}\label{interp}
		\begin{cases}
			u^-_{i+\frac12} =\omega^{\gamma}_0 \Big( \frac 1 {\gamma_0}p_0(x_{i+\frac12})  -
			\frac {\gamma_{1}} {\gamma_0} p_{1}(x_{i+\frac12})
			-\frac {\gamma_{2}} {\gamma_0} p_{2}(x_{i+\frac12}) \Big)
			+\omega^{\gamma}_{1} p_{1}(x_{i+\frac12})+\omega^{\gamma}_{2} p_{2}(x_{i+\frac12}),\\
			v^-_{i+\frac12} =p'_0(x_{i+\frac12}),
		\end{cases}
	\end{align}
	and the modified derivative $\tilde{v}_i$ is obtained by
	\begin{equation}\label{limiter}
		\tilde{v}_i = \omega^{d}_0 \Big(
		\frac{1}{d_0}p'_0(x_i) - \frac{d_1}{d_0} p'_{1}(x_i)-\frac{d_2}{d_0} p'_{2}(x_i) \Big)
		+ \omega^{d}_1 p'_1(x_i)+ \omega^{d}_2 p'_2(x_i),
	\end{equation}
	where $\omega^{d}_{\ell},~\ell=0,1,2$ are the nonlinear weights  corresponding to the linear weights $d_\ell,~\ell=0,1,2$ with the formula \eqref{nonlinear}.
{
The linear weights $\{d_0,~d_1,~d_2\}$ in HWENO limiter \eqref{limiter} and $\{\gamma_0,~\gamma_1,~\gamma_2\}$ in HWENO interpolation \eqref{interp} can be the same positive numbers in general.
Two sets of linear weights are used here to make the linear weights in HWENO reconstruction and HWENO limiter independent of each other,
thereby giving  the scheme more flexibility.
The performance of the HWENO-I scheme with different linear weights in HWENO reconstruction and HWENO limiter, respectively, is discuss in \S\ref{sec:numerical} (cf. Example \ref{blast-1d}).
}

It is worth mentioning that \eqref{limiter} is the same as the HWENO limiter in \cite{HWENO-R,HWENO-L}.
	However, compared with the HWENO-L and HWENO-R schemes \cite{HWENO-R,HWENO-L}, the modification of the derivative value in our scheme uses the same candidate stencils {(cf. Fig \ref{HWENOI-s} and Fig. \ref{HWENOR-s})} and Hermite interpolation polynomials of the reconstruction $u^-_{i+\frac12}$ and $v^-_{i+\frac12}$, and can avoid the computation of the smoothness indicator $\beta_\ell,~\ell=1,2,3$ (which have been obtained in the reconstruction of $u^{-}_{i+\frac12}$).

	\subsection{Central difference approximation for $\mathcal{D}_{i+\frac12}(f)$ and $\mathcal{D}_{i+\frac12}(h)$}
	\label{sec:Dfapp}
	
	The previous interpolation-based HWENO scheme of Liu and Qiu \cite{LiuQiu2014JSC-2} splits the fluxes $f$ and $h$ into two parts for the stability and upwind performance of their scheme.
	Unlike \cite{LiuQiu2014JSC-2}, our proposed HWENO-I scheme does not need the flux splitting and directly uses the high-order central difference approximation of $\mathcal{D}_{i+\frac12}(f)$ and $\mathcal{D}_{i+\frac12}(h)$.
	
	Consider the Hermite interpolation based on values $\{f_{i-1}$,$f_i$, $f_{i+1}$, $f_{i+2}$; $h_{i-1}$, $h_{i+2}\}$, and there is a unique fifth-degree polynomial $P(x)$ such that,
	\begin{equation}\label{Pc}
		P(x): ~\begin{cases}
			P(x_{i+\ell})=f_{i+\ell},  \quad \ell=-1,0,1,2,\\
			P'(x_{i+\ell})= h_{i+\ell},\quad\ell=-1,2.\\
		\end{cases}
	\end{equation}
	Then, we have
	\begin{equation}
		\begin{split}\label{Cd-1}
			&\Delta x^{2}P''~(x_{i+\frac12})
			= \frac{9}{8}\left( f_{i-1}-f_i-f_{i+1}+f_{i+2}\right)+\frac{5}{12}\Delta x \left(h_{i-1}-h_{i+2}\right),
			\\&
			\Delta x^{2}P^{(3)}(x_{i+\frac12})
			=-\frac{1}{\Delta x}\left( \frac{101}{36}f_{i-1}-\frac{27}{4}f_i+\frac{27}{4}f_{i+1}-\frac{101}{36}f_{i+2}\right)
			-\frac{5}{6} \left(h_{i-1}+h_{i+2}\right),
			\\&
			\Delta x^{4}P^{(4)}(x_{i+\frac12})
			=-3\left( f_{i-1}-f_i-f_{i+1}+f_{i+2}\right)-2\Delta x \left(h_{i-1}-h_{i+2}\right),
			\\&
			\Delta x^{4}P^{(5)}(x_{i+\frac12})
			=\frac{1}{\Delta x}\left( \frac{130}{9}f_{i-1}-30f_i+30f_{i+1}-\frac{130}{9}f_{i+2}\right)
			+\frac{20}{3} \left(h_{i-1}+h_{i+2}\right).
		\end{split}
	\end{equation}
	We define
	\begin{equation}\label{DfDh}
		\begin{cases}
			\mathcal{D}_{i+\frac12}(f)=-\frac{1}{24}\Delta x^{2}P''~(x_{i+\frac12}) +\frac{7}{5760}\Delta x^{4}P^{(4)}(x_{i+\frac12}) ,
			\\
			\mathcal{D}_{i+\frac12}(h)
			=
			-\frac{1}{24}\Delta x^{2}P^{(3)}~(x_{i+\frac12}) +\frac{7}{5760}\Delta x^{4}P^{(5)}(x_{i+\frac12}).
		\end{cases}
	\end{equation}
	Obviously, we have the following accuracy
	\begin{equation*}\label{Cd-Pd1}
		\begin{split}
			&\Delta x^{2}P''~(x_{i+\frac12}) =  \Delta x^{2}\left(\frac{\partial^{2} f}{\partial x^{2}}\right)_{i+\frac12} +\mathcal{O}(\Delta x^6),
			\\&
			\Delta x^{2}P^{(3)}(x_{i+\frac12}) =\Delta x^{2}\left(\frac{\partial^{2} h}{\partial x^{2}}\right)_{i+\frac12} + \mathcal{O}(\Delta x^5),
			\\&
			\Delta x^{4}P^{(4)}(x_{i+\frac12}) = \Delta x^{4}\left(\frac{\partial^{4} f}{\partial x^{4}}\right)_{i+\frac12} + \mathcal{O}(\Delta x^6),
			\\&
			\Delta x^{4}P^{(5)}(x_{i+\frac12}) =\Delta x^{4}\left(\frac{\partial^{4} h}{\partial x^{4}}\right)_{i+\frac12} + \mathcal{O}(\Delta x^5).
		\end{split}
	\end{equation*}
	Therefore, we obtain
	\begin{equation}\label{D6}
		\begin{cases}
			\mathcal{D}_{i+\frac12}(f) =
			-\frac{1}{24}\Delta x^{2}\left(\frac{\partial^{2} f}{\partial x^{2}}\right)_{i+\frac12}
			+\frac{7}{5760}\Delta x^{4}\left(\frac{\partial^{4} f}{\partial x^{4}}\right)_{i+\frac12}
			+\mathcal{O}(\Delta x^{6}),
			\\
			\mathcal{D}_{i+\frac12}(h) =
			-\frac{1}{24}\Delta x^{2}\left(\frac{\partial^{2} h}{\partial x^{2}}\right)_{i+\frac12}
			+\frac{7}{5760}\Delta x^{4}\left(\frac{\partial^{4} h}{\partial x^{4}}\right)_{i+\frac12}
			+\mathcal{O}(\Delta x^{5}).
		\end{cases}
	\end{equation}
	
	\begin{rem}
		The central difference approximation for $\mathcal{D}_{i+\frac12}(f)$ and $\mathcal{D}_{i+\frac12}(h)$ also can be
		obtained by using the Hermite interpolation based on values $\{f_{i-1}$, $f_i$, $f_{i+1}$, $f_{i+2}$; $h_{i}$, $h_{i+1}\}$ in general, i.e.,
		\begin{equation}\label{Pc2}
			P(x): ~\begin{cases}
				P(x_{i+\ell})=f_{i+\ell},  \quad \ell=-1,0,1,2,\\
				P'(x_{i+\ell})= h_{i+\ell},\quad\ell=0,1,\\
			\end{cases}
		\end{equation}
		and
		\begin{equation}\label{Cd-2}
			\begin{split}
				&\Delta x^{2}P''~(x_{i+\frac12})
				= -\frac{1}{8}\left( f_{i-1}-f_i-f_{i+1}+f_{i+2}\right)-\frac{5}{4}\Delta x \left(h_{i}-h_{i+1}\right),
				\\&\Delta x^{2}P^{(3)}(x_{i+\frac12})
				=\frac{1}{4\Delta x}\left( f_{i-1}+57f_i-57f_{i+1}-f_{i+2}\right)
				+\frac{15}{2} \left(h_{i}+h_{i+1}\right),
				\\& \Delta x^{4}P^{(4)}(x_{i+\frac12})
				=3\left( f_{i-1}-f_i-f_{i+1}+f_{i+2}\right)+6\Delta x \left(h_{i}-h_{i+1}\right),
				\\&
				\Delta x^{4}P^{(5)}(x_{i+\frac12})
				=-\frac{10}{\Delta x}\left( f_{i-1}+9f_i-9f_{i+1}-f_{i+2}\right)
				-60 \left(h_{i}+h_{i+1}\right).
			\end{split}
		\end{equation}
	\end{rem}
	
	\begin{rem}
		The above approximation \eqref{Cd-2} involves $h_i = h(u_i,v_i)=f'(u_i)v_i$ explicitly and $v_i$ would be extremely large near the discontinuities of the solution.
		Since the modified derivative $\tilde{v}_i$ is only used for time discretization in the fully-discrete HWENO-I scheme \eqref{fully-hweno} to keep the compactness of our HWENO scheme,
		it is a natural idea that we should use the derivative value carefully.
		Thus, we prefer to use approximation \eqref{Cd-1} for fluxes reconstruction in our computation.
	\end{rem}
	
	\vspace{10pt}
	
	To end this section, we summarize the procedure of the 1D HWENO-I scheme. Given the solutions $\{u_{i}^{n},v_{i}^{n}\}$ at time $t_n$, we can obtain $\{u_{i}^{n+1},v_{i}^{n+1}\}$ at time $t_{n+1}$ by solving \eqref{fully-hweno}.
	For each stage of the Runge-Kutta scheme:
	\begin{itemize}
		\item[]{\bf Step 1.} Reconstruct \{$u^{\pm}_{i+\frac{1}{2}},~v^{\pm}_{i+\frac{1}{2}}\} $ and the modified derivative $\tilde{v}_i$, see in \S\ref{sec:HWENO-I-1d}.
		\item[]{\bf Step 2.} Compute $\{ \mathcal{D}_{i+\frac12}(f), ~ \mathcal{D}_{i+\frac12}(h) \}$ by \eqref{Cd-1}-\eqref{DfDh}.
		\item[]{\bf Step 3.} Compute the numerical fluxes  $\{\hat f_{i+\frac{1}{2}},~\hat h_{i+\frac{1}{2}} \}$ by \eqref{flux-1d}.
		\item[]{\bf Step 4}.  Compute  $\mathcal{L}^1_{i}(u,v) = - \frac 1 {\Delta x} \big( \hat f_{i+\frac12}-\hat f_{i-\frac12}\big)$ and $\mathcal{L}^2_{i}(u,v) = - \frac 1 {\Delta x} \big( \hat h_{i+\frac12}-\hat h_{i-\frac12}\big)$.
	\end{itemize}
	
	\begin{rem}
		For the system case, such as the one-dimensional compressible Euler equations, the HWENO interpolation is implemented based on the local characteristic decomposition \cite{WENOJS1996}.
	\end{rem}
	
	\begin{rem}\label{BCs1d}
		We address the treatment of boundary conditions for solution $U=(\rho,~\rho \mu,~ E)$ and its derivative variables $V=U_x$ in one-dimensional compressible Euler equations \eqref{euler1d}.
		Denote the solution at ghost cell as $\bar{U}=(\bar{\rho},~\bar{\rho}\bar{ \mu}, ~\bar{E})$ and derivative variable as $\bar{V}$. For the  reflective boundary condition, we have
		\begin{equation}
\bar{\rho}=\rho,\quad \bar{\mu}=-\mu, \quad ~\bar{P} =P,\quad \bar{x} = -x,
		\end{equation}
		and
		\begin{equation}
			\frac{\partial \bar{\rho}}{\partial \bar{x}}
			=\frac{\partial  \rho}{\partial \bar{x}}
			= -\frac{\partial  \rho}{\partial x},
			\qquad
			\frac{\partial \bar{\mu}}{\partial \bar{x}}=-\frac{\partial  \mu}{\partial \bar{x}}=\frac{\partial  \mu}{\partial x},
			\qquad
			\frac{\partial \bar{P}}{\partial \bar{x}}=
			\frac{\partial P}{\partial \bar{x}}=-\frac{\partial  P}{\partial x}
		\end{equation}
		Thus, we have
		\begin{equation}\label{RBC-1d}
			\begin{split}
				&\bar{U}:~~\bar{\rho}=\rho,\qquad \bar{\rho}\bar{\mu}=-\rho\mu,  \qquad ~\bar{E} =E,\\
				&\bar{V}:~~\frac{\partial \bar{\rho}}{\partial \bar{x}}
				= -\frac{\partial  \rho}{\partial x},
				\qquad
				\frac{\partial \bar{\rho}\bar{\mu}}{\partial \bar{x}}=\frac{\partial \rho\mu}{\partial x},
				\qquad
				\frac{\partial \bar{E}}{\partial \bar{x}}=-\frac{\partial  E}{\partial x}.
			\end{split}
		\end{equation}
	\end{rem}
	
	\section{Interpolation-based HWENO scheme in 2D}
	\label{sec:HWENO2d}
	
	In this section, we extend the proposed HWENO-I scheme to solve the two-dimensional hyperbolic conservation laws. We first consider the two-dimensional scalar hyperbolic conservation laws
	\begin{equation}\label{EQ2}
		\begin{cases}
			u_t+ f(u)_x+g(u)_y=0, \\
			u(x,y,0)=u_0(x,y).
		\end{cases}
	\end{equation}
	Introduce the new variable $v = u_x$ and $w=u_y$, {and take partial derivative with respect to $x$ and $y$, respectively, on both sides of \eqref{EQ2}. Then, we can get an artificial first-order system as}
	\begin{equation}\label{EQ3}
		\begin{cases}
			u_t+ f(u)_x+g(u)_y=0, \quad u(x,y,0)=u_0(x,y),\\
			v_t+h(u,v)_x+\xi(u,v)_y=0, \quad v(x,y,0)=v_0(x,y),\\
			w_t+\eta(u,w)_x+\theta(u,w)_y=0, \quad w(x,y,0)=w_0(x,y),
		\end{cases}
	\end{equation}
	where
	\begin{align*}
		&h(u,v)=f'(u)u_x =f'(u)v, \quad \xi(u,v)=g'(u)u_x = g'(u)v,
		\\&\eta(u,w)=f'(u)u_y =f'(u)w, \quad \theta(u,w)=g'(u)u_y =g'(u)w.
	\end{align*}
	
	Let $I_{i} = [x_{i-\frac12},x_{i+\frac12}]$, $J_{j}=[y_{j-\frac12},y_{j+\frac12}]$, $I_{i,j} = I_{i}\times J_j $, $\Delta x=x_{i+\frac12}-x_{i-\frac12}$, $\Delta y=y_{j+\frac12}-y_{j-\frac12}$, and $(x_i,y_j)$ is the center of the element $I_{i,j}$.
	The semi-discrete high-order conservative finite difference scheme of \eqref{EQ3} is
	\begin{equation}
		\label{ode2H}
		\begin{cases}
			\frac{d}{dt}u_{i,j}(t) =\mathcal{L}^1_{i,j}(u,v,w) :=- \frac 1 {\Delta x} \Big( \hat f_{i+\frac12,j}-\hat f_{i-\frac12,j}\Big)
			-\frac1{\Delta y}\Big( \hat g_{i,j+\frac12}-\hat g_{i,j-\frac12}\Big),\\
			\frac{d}{dt}v_{i,j}(t) =\mathcal{L}^2_{i,j}(u,v,w) :=- \frac 1 {\Delta x} \Big( \hat h_{i+\frac12,j}-\hat h_{i-\frac12,j}\Big)
			-\frac1{\Delta y}\Big( \hat \xi_{i,j+\frac12}-\hat \xi_{i,j-\frac12}\Big),\\
			\frac{d}{dt}w_{i,j}(t) =\mathcal{L}^3_{i,j}(u,v,w) :=- \frac 1 {\Delta x} \Big( \hat \eta_{i+\frac12,j}-\hat \eta_{i-\frac12,j}\Big)
			-\frac1{\Delta y}\Big( \hat \theta_{i,j+\frac12}-\hat \theta_{i,j-\frac12}\Big).
		\end{cases}
	\end{equation}
	Here, the numerical fluxes $ \hat f_{i\pm\frac12,j}$, $\hat h_{i\pm\frac12,j}$, $\hat g_{i,j\pm\frac12}$ and $\hat \theta_{i,j\pm\frac12}$ are reconstructed with a dimension-by-dimension manner by the same procedure in \S\ref{sec:HWENO}.
	For the mixed derivative terms  $\hat \eta_{i\pm\frac12,j}$ and $\hat \xi_{i,j\pm\frac12}$, we adopt the high-order linear approximation \cite{HWENO-R,HWENO-L} directly without the flux splitting
	\begin{equation}\label{mixed}
		\begin{split}
			&\hat \eta_{i+\frac12,j}=
			-\frac{1}{12}\eta_{i-1,j}+\frac{7}{12}\eta_{i,j}
			+\frac{7}{12}\eta_{i+1,j}-\frac{1}{12}\eta_{i+2,j},\\
			&\hat \xi_{i,j+\frac12}=
			-\frac{1}{12}\xi_{i,j-1}+\frac{7}{12}\xi_{i,j}
			+\frac{7}{12}\xi_{i,j+1}-\frac{1}{12}\xi_{i,j+2}.
		\end{split}
	\end{equation}
	
	Based on the third-order SSP Runge-Kutta scheme, we have the fully-discrete scheme as, for $i=1,...,N_x,~j=1,...,N_y$,
	\begin{equation}\label{2Dfully-hweno}
		\begin{cases}
			\begin{cases}
				u^{(1)}_{i,j} =u^n_{i,j}+ \Delta t \mathcal{L}^1_{i,j}(u^{n},v^{n},w^{n}),\\
				v^{(1)}_{i,j} =\tilde{v}_{i,j}~+ \Delta t \mathcal{L}^2_{i,j}(u^{n},v^{n},w^{n}),\\
				w^{(1)}_{i,j} =\tilde{w}_{i,j} ~+ \Delta t \mathcal{L}^3_{i,j}(u^{n},v^{n},w^{n}),
			\end{cases}\\
			\begin{cases}
				u^{(2)}_{i,j} =\frac{3}{4}u^n_{i,j}+
				\frac{1}{4}\big( u^{(1)}_{i,j}+ \Delta t \mathcal{L}^1_{i,j}(u^{(1)},v^{(1)},w^{(1)})\big),\\
				v^{(2)}_{i,j} =\frac{3}{4}\tilde{v}_{i,j} ~+
				\frac{1}{4}\big( \tilde{v}^{(1)}_{i,j} ~+ \Delta t  \mathcal{L}^2_{i,j}(u^{(1)},v^{(1)},w^{(1)}\big),\\
				w^{(2)}_{i,j} =\frac{3}{4}\tilde{w}_{i,j} ~+
				\frac{1}{4}\big( \tilde{w}^{(1)}_{i,j} ~+ \Delta t  \mathcal{L}^3_{i,j}(u^{(1)},v^{(1)},w^{(1)})\big),
			\end{cases}\\
			\begin{cases}
				u^{n+1}_{i,j} =\frac{1}{3}u^n_{i,j}+
				\frac{2}{3}\big( u^{(2)}_{i,j}+ \Delta t \mathcal{L}^1_{i,j}(u^{(2)},v^{(2)},w^{(2)})\big),\\
				v^{n+1}_{i,j} =\frac{1}{3}\tilde{v}_{i,j} ~+
				\frac{2}{3}\big( \tilde{v}^{(2)}_{i,j}~ + \Delta t  \mathcal{L}^2_{i,j}(u^{(2)},v^{(2)},w^{(2)})\big),\\
				w^{n+1}_{i,j} =\frac{1}{3}\tilde{w}_{i,j} ~+
				\frac{2}{3}\big( \tilde{w}^{(2)}_{i,j}~ + \Delta t  \mathcal{L}^3_{i,j}(u^{(2)},v^{(2)},w^{(2)})\big).
			\end{cases}
		\end{cases}
	\end{equation}

	To conclude this section, we summarize the procedure of the 2D HWENO-I scheme. Given the solutions $\{u_{i,j}^{n},v_{i,j}^{n},w_{i,j}^{n}\}$ at time $t_n$, we solve \eqref{2Dfully-hweno} to obtain $\{u_{i,j}^{n+1},v_{i,j}^{n+1},w_{i,j}^{n+1}\}$ at time $t_{n+1}$.
	For each stage of the Runge-Kutta scheme:
	\begin{itemize}
		\item[]{\bf Step 1.}  Perform the HWENO-I procedure along the $x$ direction:
		\begin{itemize}
			\item[](1.1) Reconstruct $\{u^{\pm}_{i+\frac{1}{2},j},v^{\pm}_{i+\frac{1}{2},j}\} $ and the modified derivative $\tilde{v}_{i,j}$, see in \S\ref{sec:HWENO-I-1d}.
			\item[](1.2) Compute $\{ \mathcal{D}_{i+\frac12,j}(f),~ \mathcal{D}_{i+\frac12,j}(h)\}$ by \eqref{Cd-1}-\eqref{DfDh}.
			\item[](1.3) Compute the numerical fluxes  $\{\hat f_{i+\frac{1}{2},j},\hat h_{i+\frac{1}{2},j}\}$ by \eqref{flux-1d}.
		\end{itemize}
		\item[]{\bf Step 2.}  Perform the HWENO-I procedure along the $y$ direction:
		\begin{itemize}
			\item[](2.1) Reconstruct \{$u^{\pm}_{i,j+\frac{1}{2}},w^{\pm}_{i,j+\frac{1}{2}}\} $ and the modified derivative $\tilde{w}_{i,j}$, see in \S\ref{sec:HWENO-I-1d}.
			\item[](2.2) Compute $\{\mathcal{D}_{i,j+\frac12}(g),
			\mathcal{D}_{i,j+\frac12}(\theta)\}$ by \eqref{Cd-1}-\eqref{DfDh}.
			\item[](2.3) Compute the numerical fluxes $\{\hat g_{i,j+\frac{1}{2}},\hat \theta_{i,j+\frac{1}{2}}\}$ by \eqref{flux-1d}.
		\end{itemize}
		\item[]{\bf Step 3.} Compute the numerical fluxes $\{\hat \eta_{i+\frac{1}{2},j}, \hat \xi_{i,j+\frac{1}{2}}\}$ for the mixed derivative terms by \eqref{mixed}.
		\item[]{\bf Step 4.} Compute
		\begin{align*}
			&\mathcal{L}^1_{i,j}(u,v,w) = - \frac 1 {\Delta x} \big( \hat f_{i+\frac12,j}-\hat f_{i-\frac12,j}\big) - \frac 1 {\Delta y} \big( \hat g_{i,j+\frac12}-\hat g_{i,j-\frac12}\big),
			\\&\mathcal{L}^2_{i,j}(u,v,w) = - \frac 1 {\Delta x} \big( \hat h_{i+\frac12,j}-\hat h_{i-\frac12,j}\big)- \frac 1 {\Delta y} \big( \hat \xi_{i,j+\frac12}-\hat \xi_{i,j-\frac12}\big),
			\\&\mathcal{L}^3_{i,j}(u,v,w) = - \frac 1 {\Delta x} \big( \hat \eta_{i+\frac12,j}-\hat \eta_{i-\frac12,j}\big)- \frac 1 {\Delta y} \big( \hat \theta_{i,j+\frac12}-\hat \theta_{i,j-\frac12}\big).
		\end{align*}
		
	\end{itemize}
	
	\begin{rem}
		The HWENO interpolation procedures in $x$ and $y$ directions are implemented in each local characteristic direction for the two-dimensional compressible Euler equations.
		The central difference approximations for the flux functions and the numerical flux for the mixed derivative terms are performed on each component straightforwardly.
	\end{rem}
	
	\begin{rem}
		We discuss the treatment of reflective boundary condition for solution $U=(\rho,~\rho \mu,~\rho\nu,~ E)$ and its derivative variables $V=U_x$ and $W=U_y$ in two-dimensional compressible Euler equations \eqref{euler2d}.
		Denote the solution at ghost cell as $\bar{U}=(\bar{\rho},~\bar{\rho}\bar{ \mu},~\bar{\rho}\bar{ \nu}, \bar{E})$ and derivative variables as $\bar{V}$ and $\bar{W}$. For the  reflective boundary condition at the left or right, we have
		\begin{equation}
			\bar{\rho}=\rho,\quad \bar{\mu}=-\mu, \quad \bar{\nu}=\nu, \quad ~\bar{P} =P,\quad \bar{x} = -x,\quad \bar{y} = y.
		\end{equation}
		Thus,
		\begin{equation}\label{RBC-2d-x}
			\begin{split}
				&\bar{U}:~~~~\bar{\rho}=\rho,
				\qquad ~~~~~~\bar{\rho}\bar{\mu}=-\rho\mu,
				\qquad ~~~~~\bar{\rho}\bar{\nu}=\rho\nu,
				\qquad ~~~~~~\bar{E} =E,\\
				&\bar{V}:~~\frac{\partial \bar{\rho}}{\partial \bar{x}}
				= -\frac{\partial  \rho}{\partial x},
				\qquad
				\frac{\partial \bar{\rho}\bar{\mu}}{\partial \bar{x}}=\frac{\partial \rho\mu}{\partial x},
				\qquad ~~
				\frac{\partial \bar{\rho}\bar{\nu}}{\partial \bar{x}}=-\frac{\partial \rho\nu}{\partial x},
				\qquad
				\frac{\partial \bar{E}}{\partial \bar{x}}=-\frac{\partial  E}{\partial x},\\
				&\bar{W}:~~\frac{\partial \bar{\rho}}{\partial \bar{y}}
				= \frac{\partial  \rho}{\partial y},
				\qquad ~~
				\frac{\partial \bar{\rho}\bar{\mu}}{\partial \bar{y}}=-\frac{\partial \rho\mu}{\partial y},
				\qquad
				\frac{\partial \bar{\rho}\bar{\nu}}{\partial \bar{y}}=\frac{\partial \rho\nu}{\partial y},
				\qquad ~~
				\frac{\partial \bar{E}}{\partial \bar{y}}=\frac{\partial  E}{\partial y}.
			\end{split}
		\end{equation}
		Similarly, for the  reflective boundary condition at the bottom or top, we have
		\begin{equation}
			\bar{\rho}=\rho,\quad \bar{\mu}=\mu,
			\quad \bar{\nu}=-\nu,
			\quad ~\bar{P} =P,\quad \bar{x} = x,\quad \bar{y} = -y,
		\end{equation}
		and
		\begin{equation}\label{RBC-2d-y}
			\begin{split}
				&\bar{U}:~~~~\bar{\rho}=\rho,
				\qquad~~~~~~~ \bar{\rho}\bar{\mu}=\rho\mu,
				\qquad~~~~~~ \bar{\rho}\bar{\nu}=-\rho\nu,
				\qquad ~~~~~\bar{E} =E,\\
				&\bar{V}:~~\frac{\partial \bar{\rho}}{\partial \bar{x}}
				= \frac{\partial  \rho}{\partial x},
				\qquad ~~
				\frac{\partial \bar{\rho}\bar{\mu}}{\partial \bar{x}}=\frac{\partial \rho\mu}{\partial x},
				\qquad ~~
				\frac{\partial \bar{\rho}\bar{\nu}}{\partial \bar{x}}=-\frac{\partial \rho\nu}{\partial x},
				\qquad
				\frac{\partial \bar{E}}{\partial \bar{x}}=\frac{\partial  E}{\partial x},\\
				&\bar{W}:~~\frac{\partial \bar{\rho}}{\partial \bar{y}}
				= -\frac{\partial  \rho}{\partial y},
				\qquad
				\frac{\partial \bar{\rho}\bar{\mu}}{\partial \bar{y}}=-\frac{\partial \rho\mu}{\partial y},
				\qquad
				\frac{\partial \bar{\rho}\bar{\nu}}{\partial \bar{y}}=\frac{\partial \rho\nu}{\partial y},
				\qquad ~~
				\frac{\partial \bar{E}}{\partial \bar{y}}=-\frac{\partial  E}{\partial y}.
			\end{split}
		\end{equation}
	\end{rem}

	\section{Numerical experiments}
	\label{sec:numerical}
	
	In this section, we present the numerical results to show the good performance of the proposed fifth-order finite difference interpolation-based Hermite WENO (HWENO-I) scheme for one- and two-dimensional hyperbolic conservation laws.
	For comparisons, we mainly consider the HWENO-R scheme \cite{HWENO-R} since this scheme was developed recently and has an excellent performance in efficiency, robustness, and resolution than other finite difference HWENO schemes, e.g., HWENO-M \cite{ HWENO-M}, HWENO-L \cite{HWENO-L}, and {HWENO-A} \cite{LiuQiu2014JSC-2}.
In the simulation, the CFL number is taken as $0.6$ for all simulations. It is worth pointing out that the CFL number is taken as $0.2$ in the {HWENO-A} scheme \cite{LiuQiu2014JSC-2} for stability.
	We take the linear weights $\{\gamma_0,~\gamma_1=\gamma_2=(1-\gamma_0)/2\}$ in the HWENO interpolation/reconstruction with $\gamma_0=0.95$
	for one dimension and $\gamma_0=0.99$ for two dimensions, unless otherwise stated.
{
The linear weights $\{d_0=0.9,~d_1=d_2=(1-d_0)/2\}$ in the HWENO limiter for both HWENO-I and HWENO-R schemes as same in \cite{HWENO-R}.}
	To show the efficiency of the scheme, the CPU time is measured in seconds on MATLAB 2023a under the environment of Macbook Air 2020, with a M1 chip and 8GB of RAM.

	\begin{example} \label{accuracy-burgers1d}
		(Accuracy test of the one-dimensional Burgers' equation)
	\end{example}
	This example is used to verify the fifth-order accuracy and efficiency of the proposed HWENO-I scheme for the one-dimensional nonlinear Burgers' equation:
	\begin{equation*}\label{1dbugers}
		u_t+\left(\frac{u^2}{2}\right)_x=0, \quad x \in [-\pi,\pi],
	\end{equation*}
	subject to the initial condition $u(x,0)=0.5+\sin(x)$ and the periodic boundary condition.
	
\begin{figure}[H]
\centering
		\subfigure[$L^1$ error]{
			\includegraphics[width=0.45\textwidth,trim=40 0 60 10,clip]{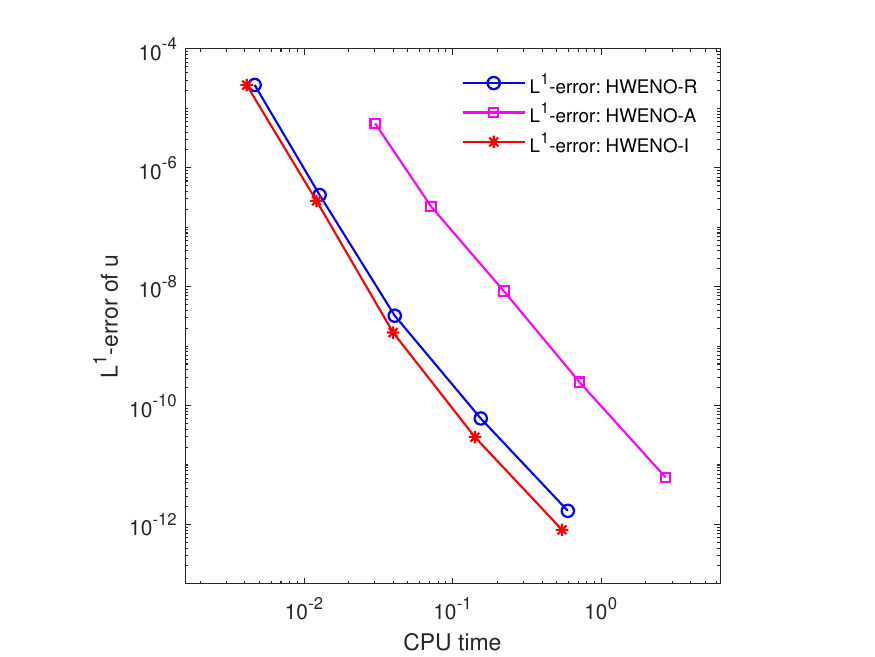}   }
		\subfigure[$L^{\infty}$ error]{
			\includegraphics[width=0.45\textwidth,trim=40 0 60 10,clip]{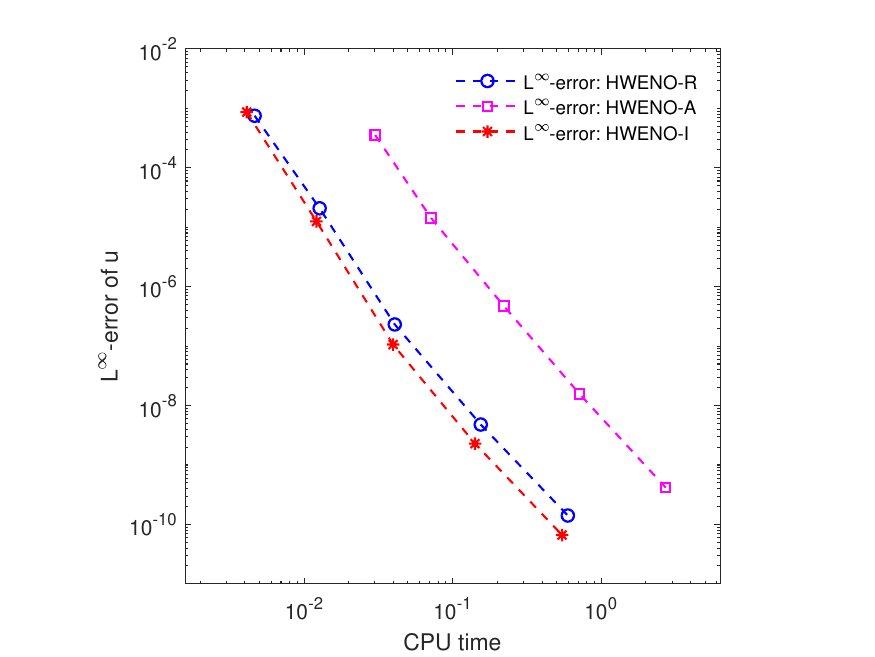}    }
		\caption{{\bf Example} \ref{accuracy-burgers1d}.
{
The errors of solution $u$ obtained by HWENO-R, HWENO-A, and HWENO-I schemes against CPU time.}
}\label{Fig:burgers1d}
\end{figure}

The final simulation time is $t=0.5$ when the solution is still smooth.
{
The errors (in $L^1$ and $L^\infty$ norm) of solution $u$ obtained by the HWENO-R and HWENO-I schemes are listed in Table~\ref{Tab:burgers1d}.
We can see that the HWENO-I scheme achieves the optimal fifth-order accuracy, and the error obtained by the HWENO-I scheme is smaller than that obtained by the HWENO-R scheme.

To show the efficiency of the scheme, the errors of solution $u$ obtained by HWENO-R, HWENO-A, and HWENO-I schemes against CPU time (measured in seconds) are plotted in Fig.~\ref{Fig:burgers1d}.
One can see that the HWENO-I scheme leads to a smaller error than the HWENO-R and HWENO-A  schemes for a fixed CPU time,
therefore the HWENO-I scheme is more accurate than the HWENO-R and HWENO-A schemes.
Furthermore, the HWENO-I scheme is more efficient than the HWENO-R scheme, and much more efficient than the HWENO-A scheme.}
	
\begin{table}[H]
\centering
\caption{{\bf Example}~\ref{accuracy-burgers1d}.
{The errors of solution $u$ obtained by the HWENO-R and HWENO-I schemes.}}
\vspace{8pt}
\begin{tabular} {c|cccc|cccc}
\toprule
& \multicolumn{4}{c}{HWENO-R}&\multicolumn{4}{c}{HWENO-I}\\
$N_x$& $L^1$ error &  order & $L^\infty$ error & order &
$L^1$ error &  order & $L^\infty$ error  &order  \\
\midrule
40	&2.480E-05	&		&7.520E-04	&	    &2.463E-05	&		&8.662E-04	&		\\
80	&3.485E-07	&6.153 	&2.085E-05	&5.173 	&2.786E-07	&6.466 	&1.255E-05	&6.109 	\\
160	&3.248E-09	&6.745 	&2.316E-07	&6.492 	&1.695E-09	&7.361 	&1.064E-07	&6.883 	\\
320	&6.067E-11	&5.742 	&4.793E-09	&5.594 	&2.935E-11	&5.851 	&2.290E-09	&5.538 	\\
640	&1.692E-12	&5.164 	&1.411E-10	&5.086 	&8.114E-13	&5.177 	&6.634E-11	&5.109 	\\
\bottomrule
\end{tabular}
\label{Tab:burgers1d}
\end{table}
	
\begin{example} \label{accuracy-Euler1d}
(Accuracy test of the one-dimensional Euler equations)
\end{example}
	This example is used to verify the fifth-order accuracy and efficiency of the proposed HWENO-I scheme for the one-dimensional Euler equations:
	\begin{equation}\label{euler1d}
		\frac{\partial}{\partial t}
		\begin{bmatrix*}[c]
			\rho \\
			\rho \mu \\
			E
		\end{bmatrix*}
		+
		\frac{\partial}{\partial x}
		\begin{bmatrix*}[c]
			\rho \mu \\
			\rho \mu^{2}+P \\
			\mu(E+P)
		\end{bmatrix*}
		=0,
	\end{equation}
	subject to the initial condition of $\rho(x,0)=1+0.2\sin(\pi x)$, $\mu(x,0)=1$, $P(x,0)=1$
	and the periodic boundary conditions.
	Here $\rho$ is the density, $\mu$ is the velocity, $E$ is the total energy, and $P$ is the pressure.  In the simulation of one-dimensional Euler equations, the equation of state (EOS) is defined as $E = P/(\gamma-1)+\rho \mu^2/2 $ with the ratio of specific heats $\gamma=1.4$ unless otherwise stated.
	
\begin{figure}[H]
\centering
\subfigure[$L^1$ error]{
\includegraphics[width=0.45\textwidth,trim=40 0 60 10,clip]{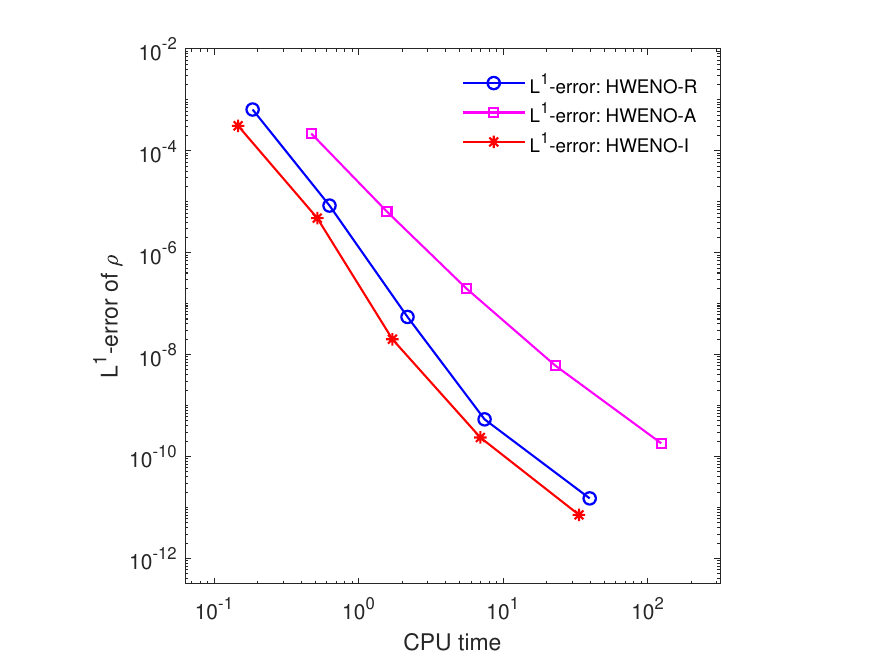}   }
\subfigure[$L^{\infty}$ error]{
\includegraphics[width=0.45\textwidth,trim=40 0 60 10,clip]{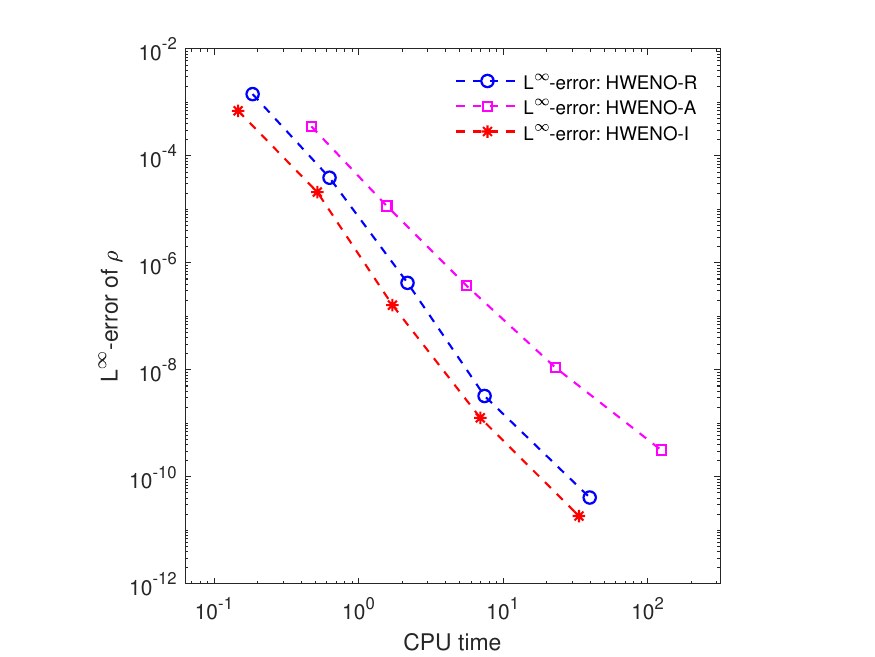}    }
\caption{{\bf Example} \ref{accuracy-Euler1d}.
{
The errors of density $\rho$ obtained by HWENO-R, HWENO-A, and HWENO-I schemes against CPU time.}}
\label{Fig:Euler1d}
\end{figure}

The computational domain is $[0, 2]$. The final simulation time is $t=2$.
{
		In Table \ref{Tab:Euler1d}, we list the errors (in $L^1$ and $L^\infty$ norm) of density $\rho$ obtained by HWENO-R and HWENO-I schemes.}
We can see that the HWENO-I scheme achieves the optimal fifth-order accuracy, and the error obtained by the HWENO-I scheme is smaller than that obtained by the HWENO-R scheme.
{
In Fig.~\ref{Fig:Euler1d}, we plot the errors of density $\rho$ obtained by HWENO-R, HWENO-A, and HWENO-I schemes against CPU time (measured in seconds).
We can see that the HWENO-I scheme has a smaller error than either the HWENO-R scheme or the HWENO-A scheme for a fixed amount of CPU time.
Thus, the HWENO-I scheme is more accurate and efficient than the HWENO-R and HWENO-A schemes in this example.}
	
\begin{table}[H]
\centering
\caption{{\bf Example}~\ref{accuracy-Euler1d}.
{The errors of density $\rho$ obtained by the HWENO-R and HWENO-I schemes.}}
\vspace{8pt}
\begin{tabular} {c|cccc|cccc}
\toprule
& \multicolumn{4}{c}{HWENO-R}&\multicolumn{4}{c}{HWENO-I}\\
$N_x$& $L^1$ error &  order & $L^\infty$ error & order &
$L^1$ error &  order & $L^\infty$ error  &order  \\
\midrule
20	&6.474E-04	&		&1.432E-03	&		&3.098E-04	&		&6.973E-04	&	 	\\
40	&8.421E-06	&6.264 	&3.914E-05	&5.194 	&4.771E-06	&6.021 	&2.104E-05	&5.050 	\\
80	&5.489E-08	&7.261 	&4.239E-07	&6.529 	&2.002E-08	&7.897 	&1.626E-07	&7.016 	\\
160	&5.359E-10	&6.678 	&3.253E-09	&7.026 	&2.363E-10	&6.404 	&1.269E-09	&7.002 	\\
320	&1.506E-11	&5.153 	&4.111E-11	&6.306 	&7.166E-12	&5.044 	&1.854E-11	&6.097 	\\
\bottomrule
\end{tabular}
\label{Tab:Euler1d}
\end{table}	
		
	\begin{example}\label{ShuOsher-1d}
		(Shu-Osher problem of the one-dimensional Euler equations)
	\end{example}
	In this example, we consider the Shu-Osher problem of the one-dimensional Euler equations \eqref{euler1d} with the initial condition:
	\begin{equation*}
		(\rho, \mu, P)=\begin{cases}
			(3.857143,~2.629369,~10.333333),&\text { for } x<-4, \\
			(1+0.2 \sin (5 x), ~0,~1), &\text { for } x \geq-4,
		\end{cases}
	\end{equation*}
which is a typical example of the shock interaction with the entropy wave problem \cite{Shu-Osher1989}.
The solution has a moving Mach $3$ shock interacting with a sine wave and contains both shock waves and complex smooth region structures. A high-order and high-resolution scheme would show its advantage when the solution contains both shocks and complex structures.

\begin{figure}[H]
       \centering
       \subfigure[density $\rho$]{
       \includegraphics[width=0.45\textwidth,trim=50 0 65 10,clip]{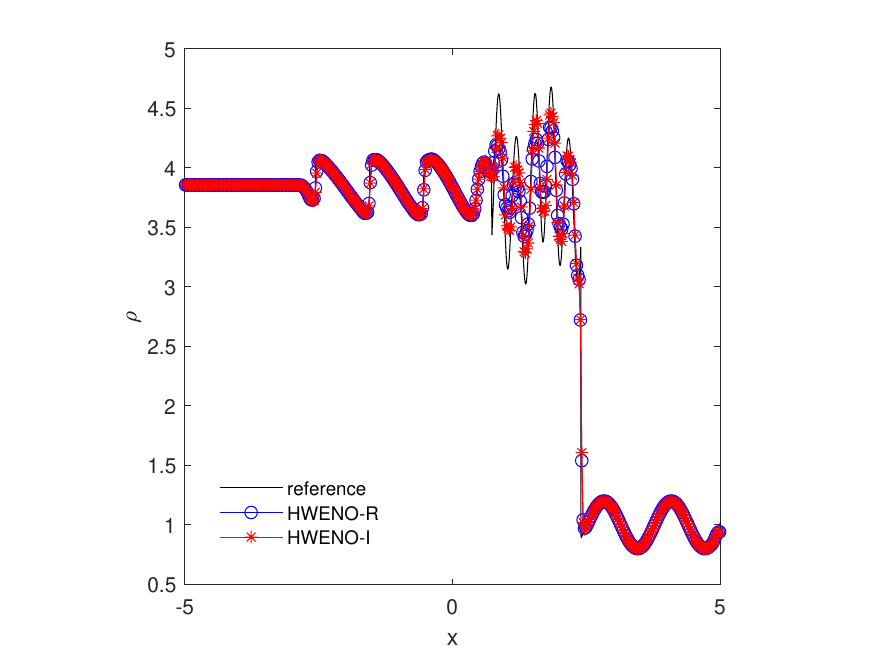}   }
      \subfigure[zoom-in (a) at $x\in(0,2.5)$]{
       \includegraphics[width=0.45\textwidth,trim=50 0 65 10,clip]{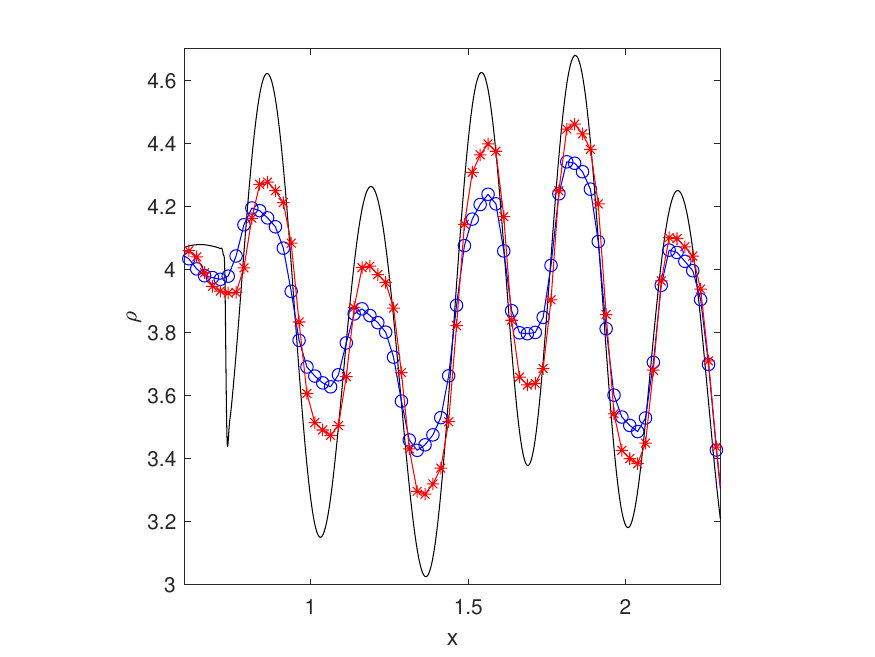}    }
 \caption{{\bf Example} \ref{ShuOsher-1d}.
 The density $\rho$ at $t=1.8$ obtained by the HWENO-R and HWENO-I schemes with $N_x=400$.}
 \label{Fig:ShuOsher1d}
 \end{figure}

The computation domain is $[-5,5]$ and the final simulation time is $t=1.8$.
The density $\rho$ at the final time obtained by the HWENO-R and HWENO-I schemes with $N_x=400$ is shown in Fig.~\ref{Fig:ShuOsher1d}.
The reference solution is obtained by HWENO-R scheme \cite{HWENO-R} with $N_x = 10000$ for comparison.
One can see that the HWENO-I scheme produces a higher resolution solution near the region $[0.8, 2.2]$ (see in Fig.~\ref{Fig:ShuOsher1d} (b)), so the HWENO-I scheme has better resolution than the HWENO-R scheme for complex smooth region structures.
	
\begin{figure}[H]
		\centering
		\subfigure[density $\rho$]{
			\includegraphics[width=0.45\textwidth]{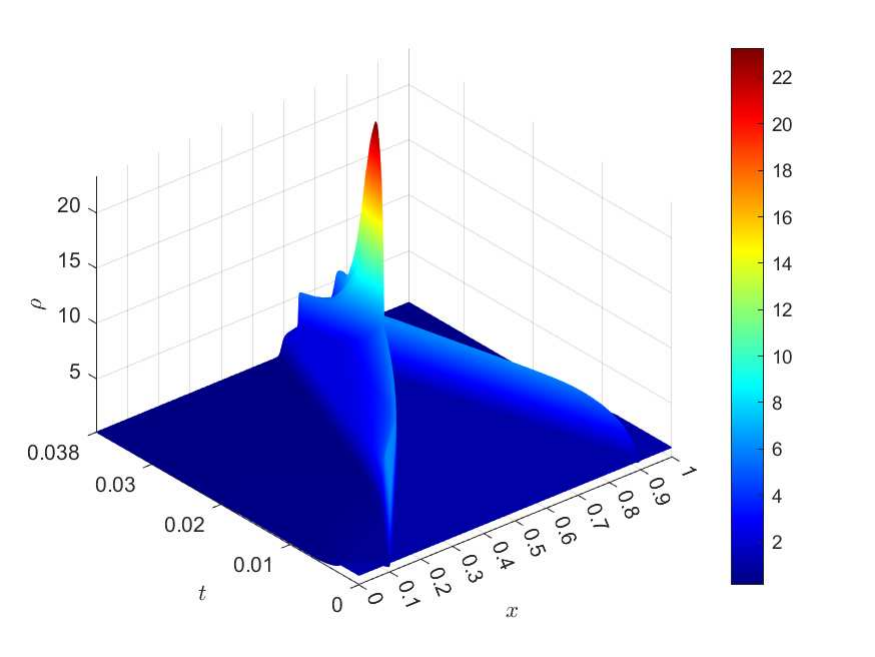}   }
		\subfigure[{$60$ contours of $\rho$:  $0.1470$ to $23.2289$}]{
			\includegraphics[width=0.45\textwidth]{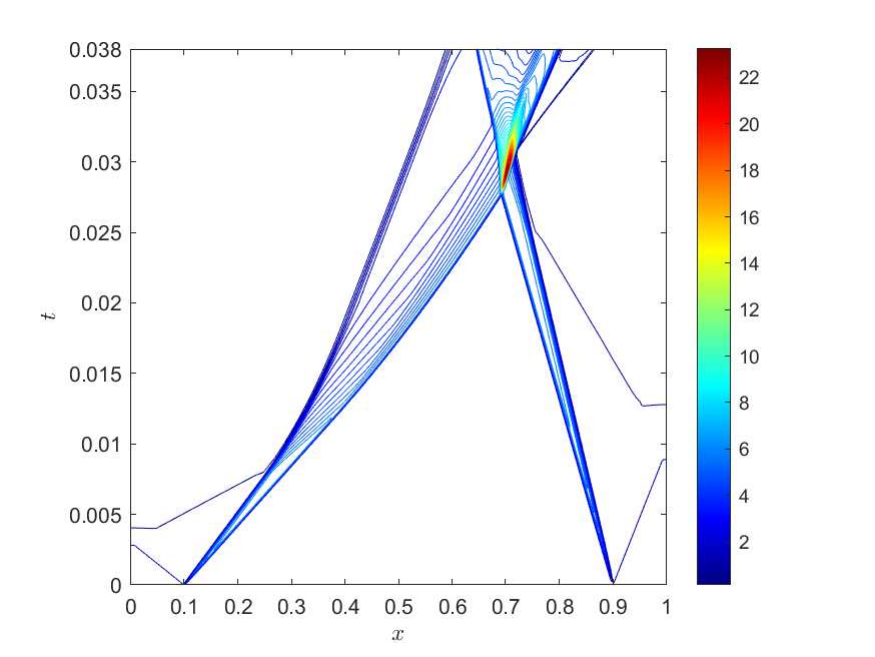}}
		\caption{{\bf Example} \ref{blast-1d}.
			The density $\rho$ on the space-time plane obtained by the HWENO-I scheme with $N_x=800$.}\label{Fig:blast1d-times}
	\end{figure}

	\begin{example}\label{blast-1d}
		(Woodward-Colella blast wave problem of the one-dimensional Euler equations)
	\end{example}
	We consider the interaction of the blast wave problem \cite{{Woodward-Colella1984}} of the one-dimensional Euler equations \eqref{euler1d} with the initial condition:
	\begin{equation*}
		(\rho, \mu, P)= \begin{cases}
			(1,~0,~1000), & \text { if } x \in[0,0.1),  \\
			(1,~0,~0.01), & \text { if } x \in[0.1,0.9), \\
			(1,~0,~100), & \text { if } x \in [0.9,1],
		\end{cases}
	\end{equation*}
	subject to the reflective boundary conditions on the left and right.  The solution involves multiple interactions of strong shocks and rarefactions with each other and with contact discontinuities. See Fig~\ref{Fig:blast1d-times}.
	
\begin{figure}[H]
		\centering
		\subfigure[density $\rho$]{
			\includegraphics[width=0.45\textwidth,trim=25 0 40 10,clip]{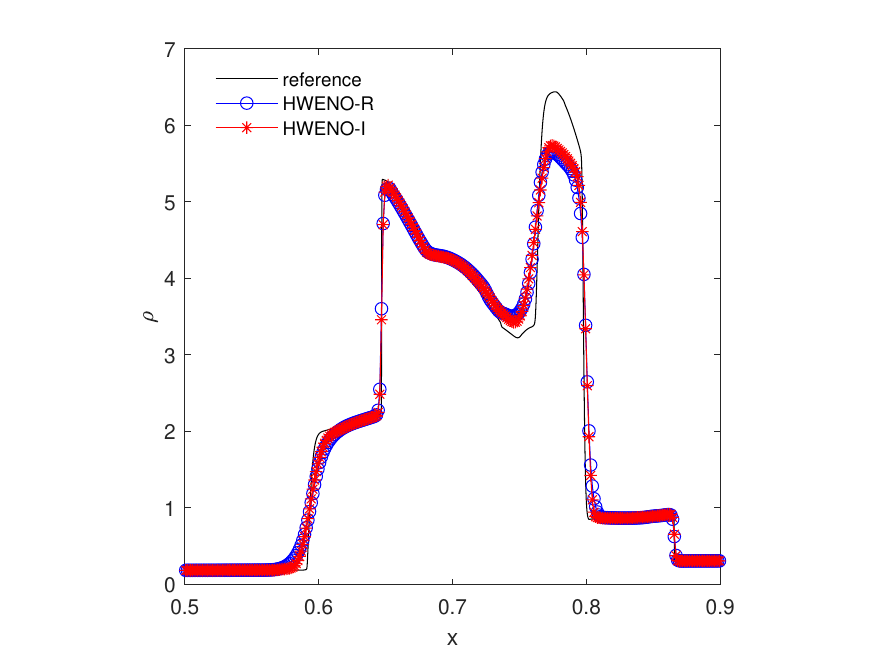}   }
		\subfigure[{zoom-in (a) at $x\in(0.55,0.63),(0.7,0.82)$}]{
			\includegraphics[width=0.45\textwidth,trim=25 0 40 10,clip]{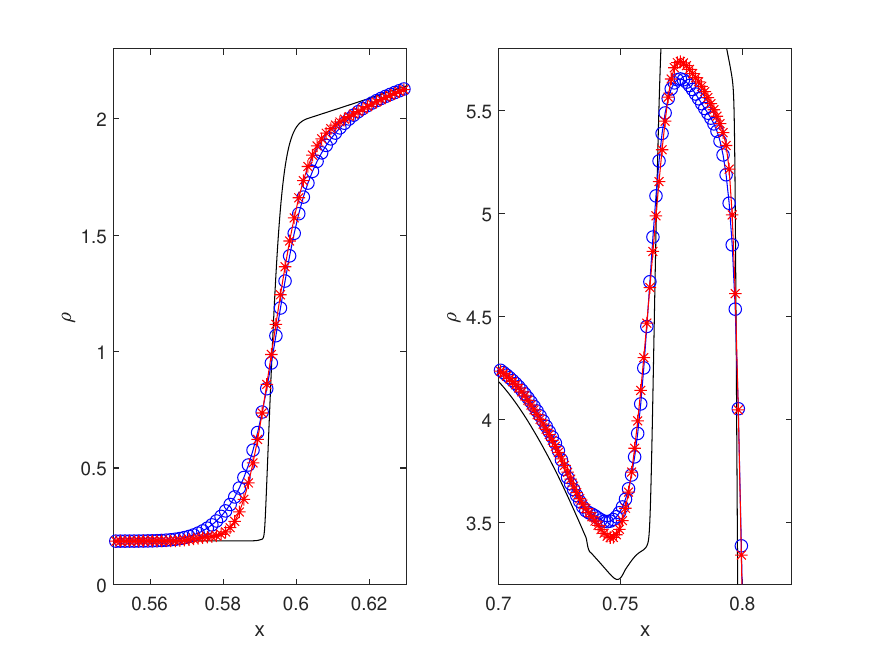}}
		\caption{{\bf Example} \ref{blast-1d}.
			The density $\rho$ at $t=0.038$ obtained by the HWENO-R and HWENO-I schemes with $N_x=800$.}\label{Fig:blast1d}
	\end{figure}

	The final time is $t=0.038$.
	The density $\rho$ obtained by the HWENO-R and HWENO-I schemes with $N_x=800$ is shown in Fig.~\ref{Fig:blast1d} in which the HWENO-R scheme computes the reference solution with $N_x=10000$ for comparison.
	From Fig.~\ref{Fig:blast1d}, we can see that the HWENO-I scheme could get sharp shock transitions in the computational fields, and the HWENO-I scheme has a higher resolution than the HWENO-R scheme near the shock and the peak.

\begin{figure}[H]
		\centering
		\subfigure[density $\rho$]{
			\includegraphics[width=0.45\textwidth,trim=50 0 65 10,clip]{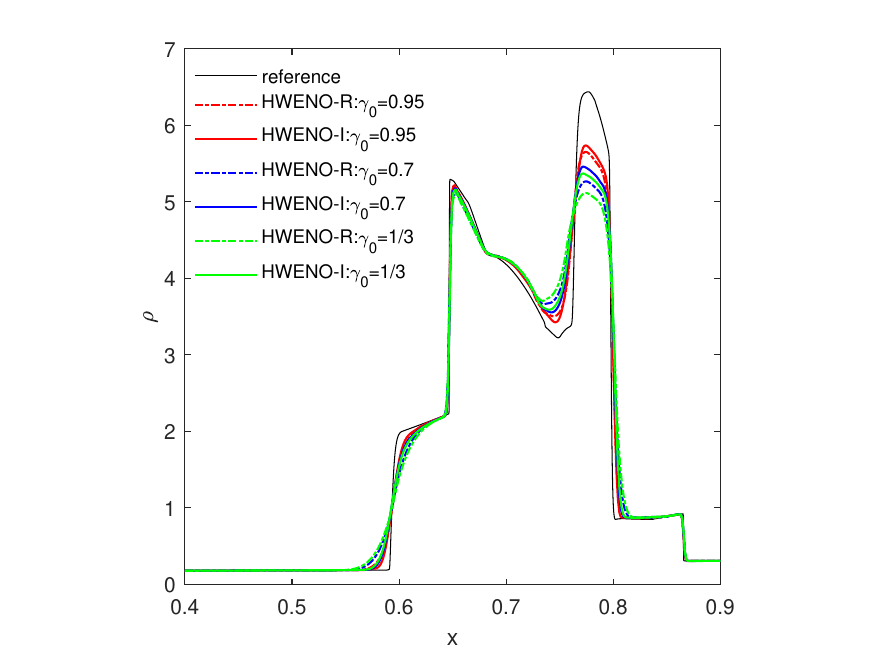}   }
		\subfigure[zoom-in (a) at $x\in(0.7, 0.8)$]{
			\includegraphics[width=0.45\textwidth,trim=50 0 65 10,clip]{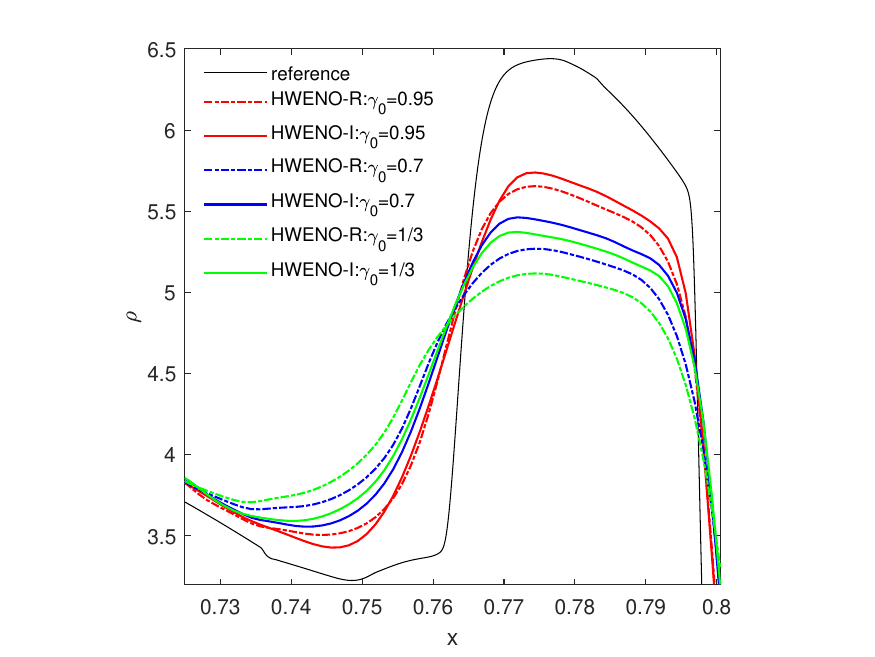}}
		\subfigure[HWENO-R: $\rho$ at $x\in(0.7, 0.8)$]{
			\includegraphics[width=0.45\textwidth,trim=50 0 65 10,clip]{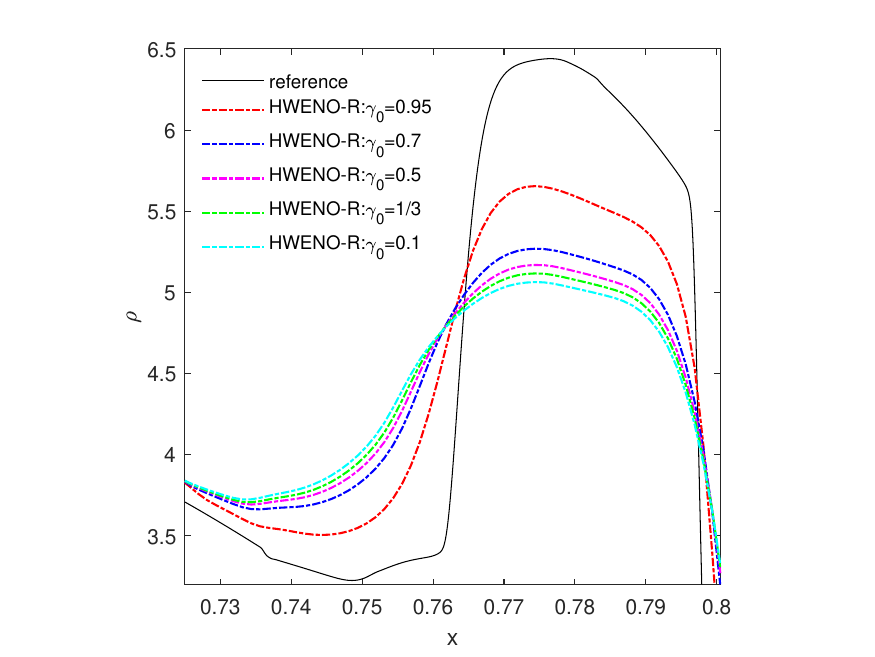}}
		\subfigure[HWENO-I: $\rho$ at $x\in(0.7, 0.8)$]{
			\includegraphics[width=0.45\textwidth,trim=50 0 65 10,clip]{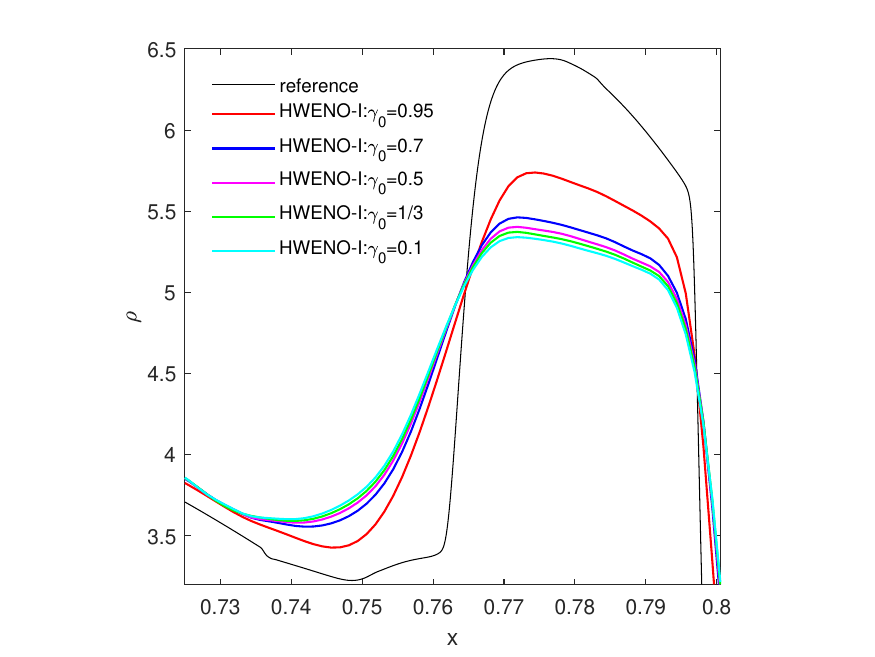}}
		\caption{{\bf Example} \ref{blast-1d}.
			The density $\rho$ (with $N_x=800$) at $t=0.038$ obtained by the HWENO-R and HWENO-I schemes with different $\gamma_0=0.95,~0.7,~0.5,~1/3,~0.1$.
The linear weights $d_0=0.9,~d_1=d_2=(1-d_0)/2\}$ are fixed in the HWENO limiter.
}
\label{Fig:blast1d-g0}
	\end{figure}

\begin{figure}[H]
		\centering
		\subfigure[$\rho$ at $x\in(0.73, 0.82)$]{
			\includegraphics[width=0.45\textwidth,trim=30 0 35 10,clip]{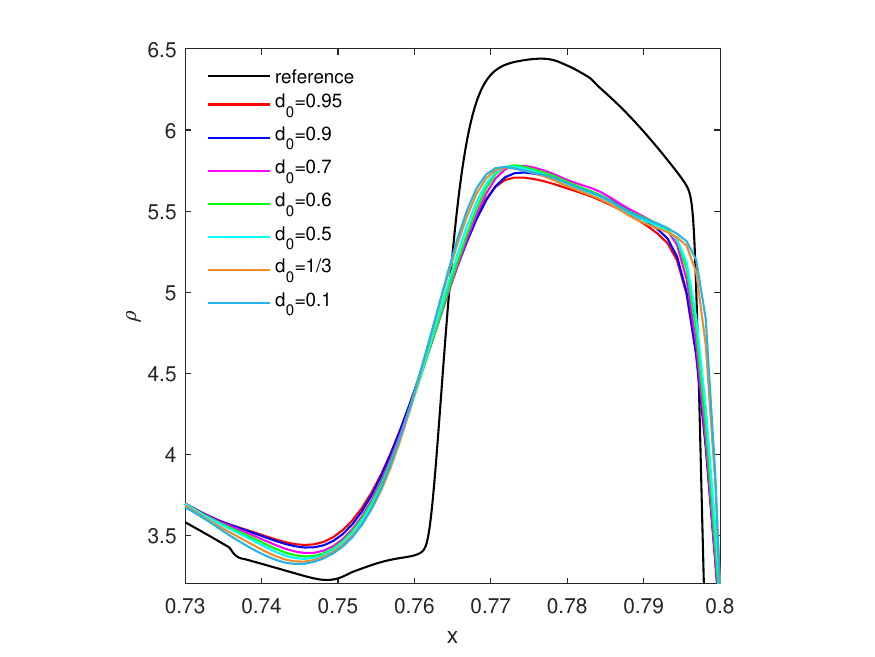}   }
		\subfigure[$\rho$ at $x\in(0.57, 0.62)\cup(0.79, 0.81)$]{
			\includegraphics[width=0.45\textwidth,trim=30 0 35 10,clip]{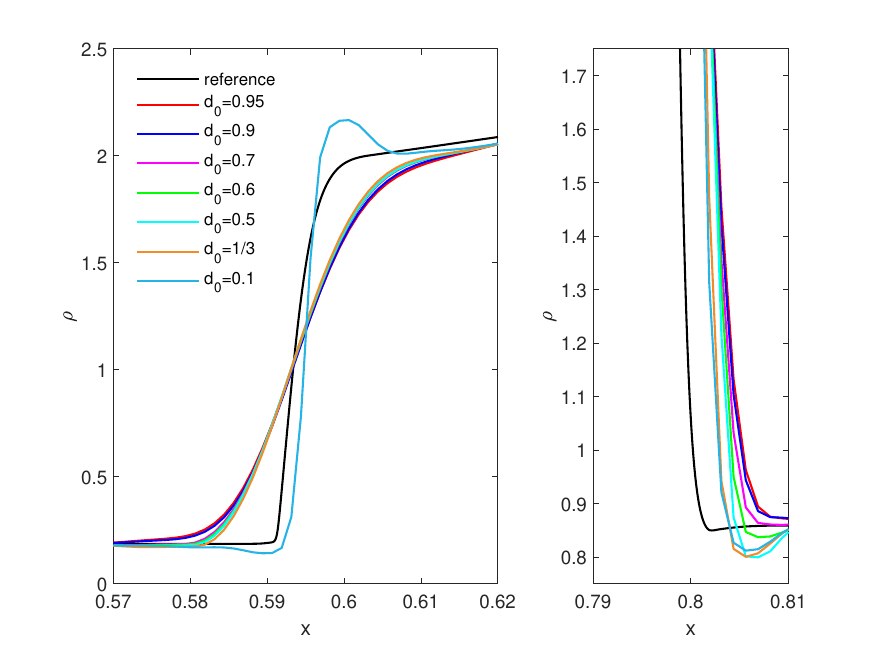}}
		\caption{{\bf Example} \ref{blast-1d}.
		{The density $\rho$ (with $N_x=800$) at $t=0.038$ obtained by the HWENO-I scheme with different $d_0=0.95,~0.9,~0.7,~0.5,~1/3,~0.1$.
The linear weights $\{\gamma_0=0.95,~\gamma_1=\gamma_2=(1-\gamma_0)/2\}$ are fixed in the HWENO interpolation/reconstruction.
}}
\label{Fig:blast1d-fixgamma}
	\end{figure}

{
Next, we study how the linear weights in HWENO reconstruction and HWENO limiter, respectively, affect the performance of the HWENO-I scheme.
First, we investigate how the linear weights $\{\gamma_0,~\gamma_1=\gamma_2=(1-\gamma_0)/2\}$ in the HWENO interpolation/reconstruction affect the performance of the HWENO-I scheme.
In Fig.~\ref{Fig:blast1d-g0},
we plot the density obtained by the HWENO-I and HWENO-R schemes with $\gamma_0=0.95,~0.7,~0.5,~1/3,~0.1$ (while the linear weights $d_0=0.9,~d_1=d_2=(1-d_0)/2\}$ in the HWENO limiter are fixed).
	From Fig.~\ref{Fig:blast1d-g0}, one can see that the resolution is decreasing as $\gamma_0$ decreases,
	and the decreasing rate of the HWENO-I scheme seems slower than that of the HWENO-R scheme from Fig.~\ref{Fig:blast1d-g0}(c) and (d).
	Moreover, we observe that the resolution of the HWENO-I scheme with $\gamma_0=1/3$ is better than that of the HWENO-R scheme with $\gamma_0=0.7$ from Fig.~\ref{Fig:blast1d-g0}(a) and (b).

Second, we investigate how the linear weights $\{d_0,~d_1=d_2=(1-d_0)/2\}$ in the HWENO limiter affect the performance of the HWENO-I scheme.
In Fig.~\ref{Fig:blast1d-fixgamma}, we show the density obtained by the HWENO-I scheme with $d_0=0.95,~0.9,~0.7,~0.5,~1/3,~0.1$ (while the linear weights $\{\gamma_0=0.95,~\gamma_1=\gamma_2=(1-\gamma_0)/2\}$ in the HWENO interpolation/reconstruction are fixed).
	From Fig.~\ref{Fig:blast1d-fixgamma}, we can obverse that the resolution of the HWENO-I scheme seems to increase slightly as $d_0$ decreases, while the robustness of the HWENO-I scheme decreases.  Especially, when $d_0$ is relatively small (e.g. $d_0< 0.7 $), the numerical oscillations appear in the solution.
To balance the resolution and robustness, we adopt the linear weights $d_0=0.9$ in the computation as in \cite{HWENO-R}.
}
	
	\begin{example} \label{accuracy-burgers2d}
		(Accuracy test of the two-dimensional Burgers' equation)
	\end{example}
	This example is used to verify the fifth-order accuracy and efficiency of the proposed HWENO-I scheme for the two-dimensional nonlinear Burgers' equation:
	\begin{equation*}\label{2dbugers}
		u_t+\left(\frac{u^2}{2}\right)_x +\left(\frac{u^2}{2}\right)_y=0, \quad (x,y) \in [-2\pi, 2\pi]\times [-2\pi, 2\pi] ,
	\end{equation*}
	subject to the initial condition $u(x,0)=0.5+\sin((x+y)/2)$ and the periodic boundary condition.
	
\begin{figure}[H]
		\centering
		\includegraphics[width=0.45\textwidth,trim=50 0 65 10,clip]
		{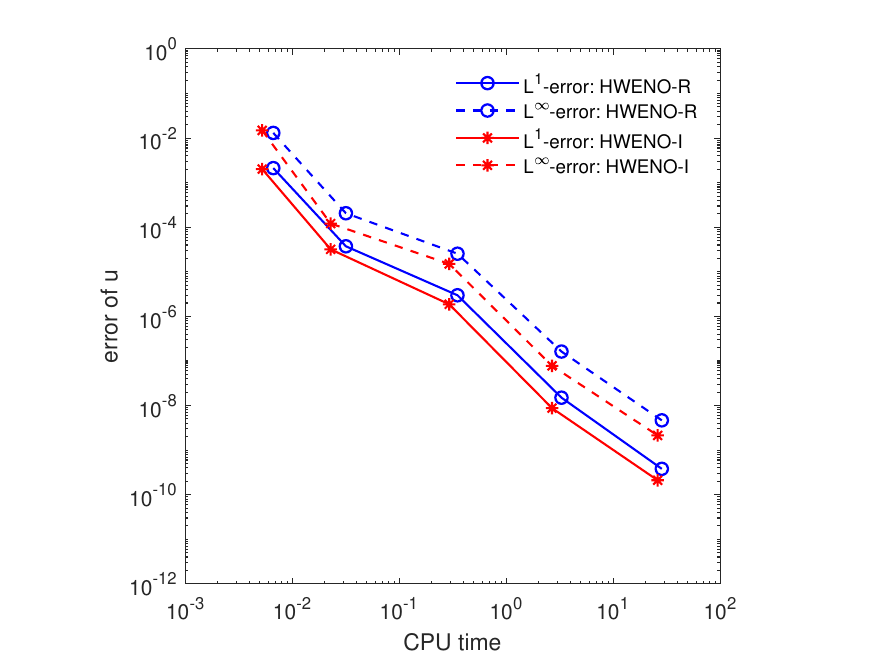}
		\caption{{\bf Example} \ref{accuracy-burgers2d}.
			{
				The errors of solution $u$ obtained by HWENO-R and HWENO-I schemes against CPU time.}}
		\label{Fig:burgers2d}
	\end{figure}

	The final simulation time is $t=0.5$ when the solution is still smooth.
	{
		In Table \ref{Tab:burgers2d}, we list the errors (in $L^1$ and $L^\infty$ norm) of solution $u$ obtained by the HWENO-R and HWENO-I schemes. The result shows that the HWENO-I achieves the fifth-order accuracy and gets smaller error than the HWENO-R scheme on the same mesh.
		In Fig. \ref{Fig:burgers2d}, we plot the errors against CPU time.
		It shows that the HWENO-I scheme is more efficient than the HWENO-R scheme in simulating this problem.}
	\begin{table}[H]
		\centering
		\caption{{\bf Example}~\ref{accuracy-burgers2d}.
			{The errors of solution $u$ obtained by the HWENO-R and HWENO-I schemes.}}
		\vspace{8pt}
		\begin{tabular} {c|cccc|cccc}
			\toprule
& \multicolumn{4}{c}{HWENO-R}&\multicolumn{4}{c}{HWENO-I}\\
$N_x=N_y$& $L^1$ error &  order & $L^\infty$ error & order &
		$L^1$ error &  order & $L^\infty$ error  &order  \\
\midrule
20  &2.122E-03&	  &1.307E-02& &2.018E-03&     &1.492E-02&	\\
40  &3.740E-05&5.826&2.062E-04&5.987&3.196E-05&5.980&1.210E-04&6.946 	\\
80  &2.980E-06&3.650&2.550E-05&3.015&1.868E-06&4.097&1.507E-05&3.006 	\\
160&1.502E-08&7.632&1.628E-07&7.291&8.754E-09&7.737&7.784E-08&7.597 	\\
320&3.805E-10&5.303&4.658E-09&5.127&2.106E-10&5.377&2.137E-09&5.187 	\\
\bottomrule
\end{tabular}
\label{Tab:burgers2d}
\end{table}
	
	\begin{example} \label{accuracy-Euler2d}
		(Accuracy test of the two-dimensional Euler equations)
	\end{example}
	This example is used to verify the fifth-order accuracy and efficiency of the proposed HWENO-I scheme for the two-dimensional Euler equations:
	\begin{equation}\label{euler2d}
		\frac{\partial}{\partial t}
		\begin{bmatrix*}[c]
			\rho \\
			\rho \mu \\
			\rho \nu \\
			E
		\end{bmatrix*}
		+
		\frac{\partial}{\partial x}
		\begin{bmatrix*}[c]
			\rho \mu \\
			\rho \mu^2 +P\\
			\rho \mu\nu \\
			\mu(E+P)
		\end{bmatrix*}
		+
		\frac{\partial}{\partial y}
		\begin{bmatrix*}[c]
			\rho \nu \\
			\rho \mu\nu \\
			\rho \nu^2+P\\
			\nu(E+P)
		\end{bmatrix*}
		=0,
	\end{equation}
	subject to the initial condition of $\rho(x,y,0)=1+0.2\sin(\pi (x+y))$, $u(x,y,0)=1$, $\nu(x,y,0)=1$, $P(x,y,0)=1$,
	and the periodic boundary conditions.
	Here $\rho$ is the density, $\mu,~\nu$ are the velocities of $x$- and $y$-direction, respectively, $E$ is the total energy and $P$ is the pressure.
	For the two-dimensional Euler equations, the EOS is defined as $E = P/(\gamma-1)+\rho (\mu^2+\nu^2)/2 $ with the ratio of specific heats $\gamma=1.4$  unless otherwise stated.
	
\begin{figure}[H]
		\centering
		\includegraphics[width=0.45\textwidth,trim=50 0 65 10,clip]
		{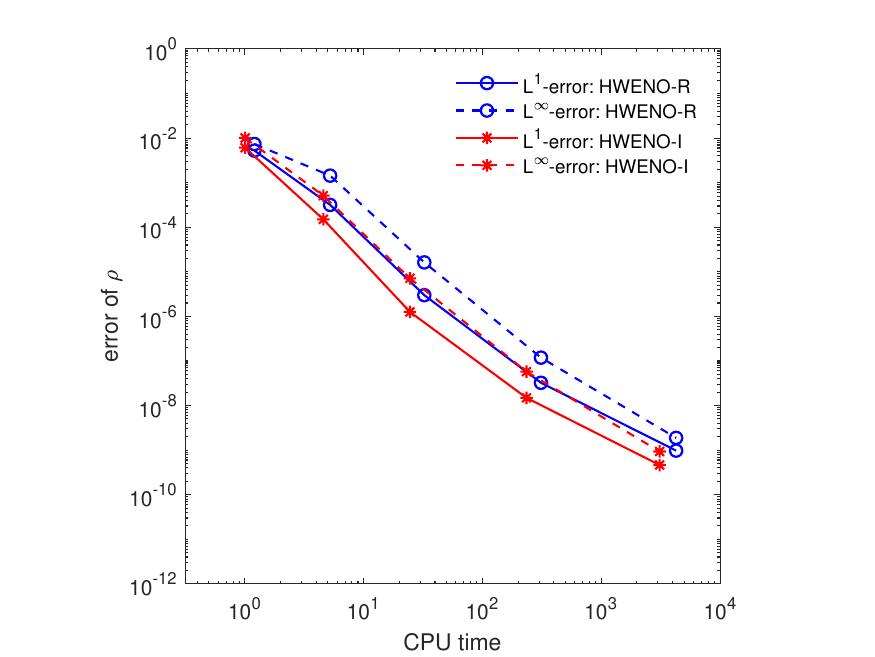}
		\caption{{\bf Example} \ref{accuracy-Euler2d}.
			{
				The errors of density $\rho$ obtained by HWENO-R and HWENO-I schemes against CPU time.}}
		\label{Fig:Euler2d}
	\end{figure}

	The computational domain is $[0, 2]\times [0,2]$.
	The final simulation time is $t=2$.
	{
		The errors (in $L^1$ and $L^\infty$ norm) of density $\rho$ obtained by HWENO-R and HWENO-I schemes are listed in Table \ref{Tab:Euler2d}, and the errors against CPU time are presented in Fig.~\ref{Fig:Euler2d}.}
One can observe that the HWENO-I scheme has the fifth-order accuracy, and is more accurate and efficient than the HWENO-R scheme in simulating the two-dimensional Euler equations.
\begin{table}[H]
\centering
\caption{{\bf Example}~\ref{accuracy-Euler2d}.
{The errors of density $\rho$ obtained by the HWENO-R and HWENO-I schemes.}}
\vspace{8pt}
\begin{tabular} {c|cccc|cccc}
\toprule
& \multicolumn{4}{c}{HWENO-R}&\multicolumn{4}{c}{HWENO-I}\\
$N_x=N_y$			& $L^1$ error &  order & $L^\infty$ error & order &
			$L^1$ error &  order & $L^\infty$ error  &order  \\
\midrule
10  &5.245E-03& 	 &7.357E-03&       &6.093E-03	&	 	&1.008E-02&	\\
20  &3.175E-04&4.046  &1.439E-03&2.354  &1.506E-04	&5.338 	&5.113E-04	&4.301 	\\
40  &2.994E-06&6.729  &1.628E-05&6.466  &1.261E-06	&6.901 	&7.161E-06	&6.158 	\\
80  &3.241E-08&6.530  &1.189E-07&7.098 &1.488E-08	&6.405 	&5.768E-08	&6.956 	\\
160 &9.750E-10&5.055  &1.882E-09&5.981 &4.630E-10	&5.006  &9.328E-10	&5.950 	\\
\bottomrule
\end{tabular}
		\label{Tab:Euler2d}
	\end{table}
	
	\begin{example}\label{doublemach}
		(Double Mach reflection problem of the two-dimensional Euler equations)
	\end{example}
	We consider the double Mach reflection problem \cite{Woodward-Colella1984} of the two-dimensional Euler equations \eqref{euler2d} over computational domain $[0,4]\times[0,1]$ which involves a pure right-moving Mach $10$ shock, initially located at $x=\frac{1}{6},~ y=0$, making a $60^{\circ}$ angle with the $x$-axis. The undisturbed air ahead of the shock has a density $\rho=1.4$ and a pressure $P=1$. Specifically, the initial condition is
	\begin{equation*}
		(\rho, \mu, \nu, P)=\begin{cases}
			\left(8, 8.25 \cos \left(\frac{\pi}{6}\right),-8.25\sin \left(\frac{\pi}{6}\right), 116.5\right),  &\text{if} ~x<\frac{1}{6}+\frac{y}{\sqrt{3}}, \\
			(1.4,~0,~0,~1) & \text { otherwise. }
		\end{cases}
	\end{equation*}
	The exact motion of the initial Mach $10$ shock is used on the top boundary.
	The exact post-shock condition and the reflective boundary are posed on bottom boundary $\left[0, \frac{1}{6}\right] \times\{0\}$ and $(\frac{1}{6},4]\times\{0\}$, respectively.
	The inflow and outflow boundary conditions are used on the left and right boundaries, respectively.
	
\begin{figure}[H]
		\centering
		\subfigure[HWENO-R]{
			\includegraphics[width=0.9\textwidth,trim=20 80 20 90,clip]{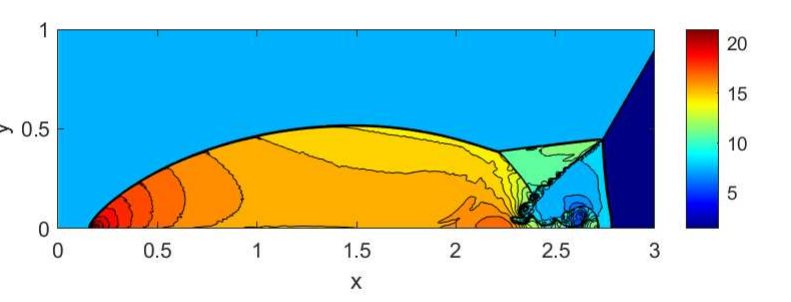}}
		\subfigure[HWENO-I]{
			\includegraphics[width=0.9\textwidth,trim=20 80 20 90,clip]{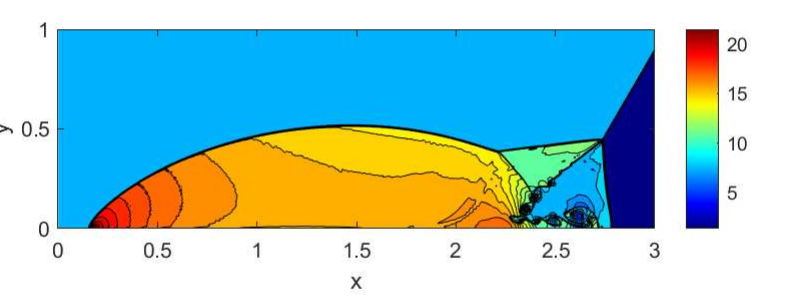}}
		\subfigure[HWENO-R at $(2.2, 2.9)\times(0,0.5)$]{
			\includegraphics[width=0.45\textwidth,trim=20 10 20 10,clip]{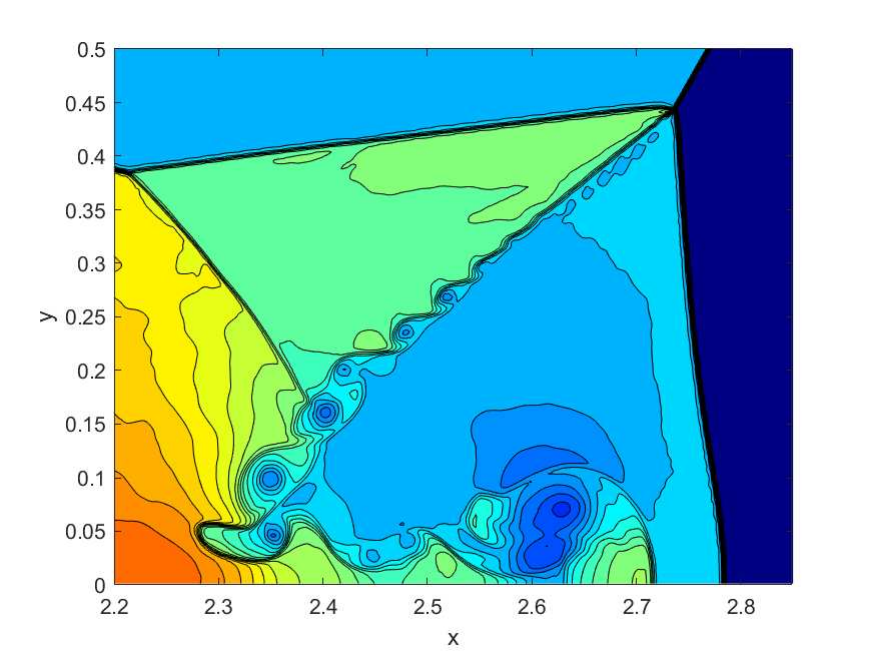}}
		\subfigure[HWENO-I at $(2.2, 2.9)\times(0,0.5)$]{
			\includegraphics[width=0.45\textwidth,trim=20 10 20 10,clip]{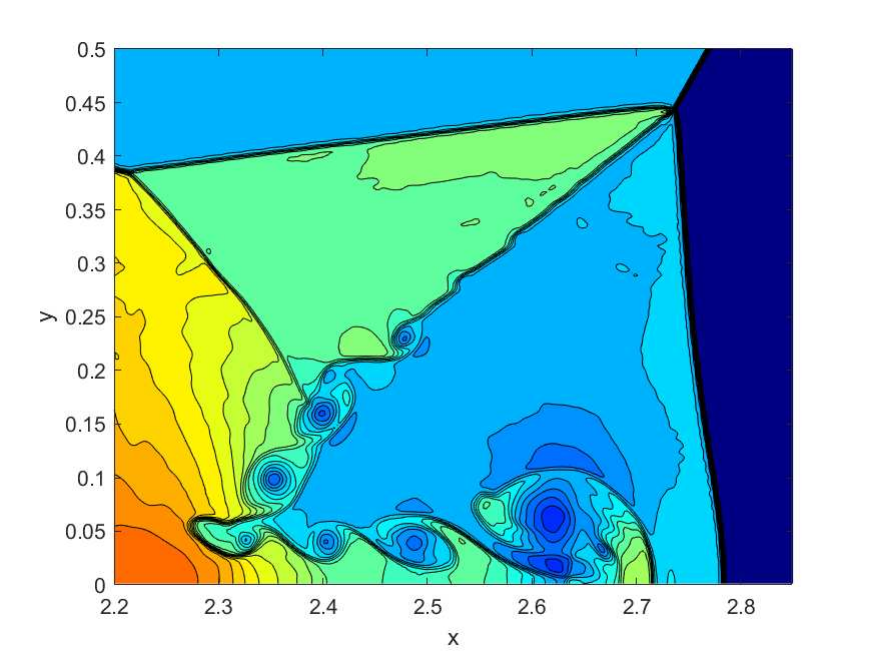}}
		\caption{{\bf Example} \ref{doublemach}.
			{
				The density $\rho$ at $t=0.2$ obtained by the HWENO-R and HWENO-I schemes with $N_x\times N_y=1600\times 400$.}
		}\label{Fig:doublemach2d}
	\end{figure}

	The density $\rho$ at $t=0.2$ obtained by the HWENO-R and HWENO-I schemes with $N_x\times N_y=1600\times 400$ is presented in Fig.\ref{Fig:doublemach2d}.
	One can see that both schemes can capture complicated structures around the double Mach stems and the HWENO-I scheme has a higher resolution than the HWENO-R scheme near the region $(2.4, 2.7)\times(0,0.1)$.
	It is worth pointing out that the total CPU time of the HWENO-R scheme with $N_x\times N_y=1600\times 400$ is $29805.41s$ while the HWENO-I scheme only costs $21683.64s$.
	The CPU time needed by the HWENO-I scheme is only about $73\%$ of that used by the HWENO-R scheme, so the HWENO-I scheme is more efficient.
	
	\begin{example}\label{forwardstep}
		(Forward step problem of the two-dimensional Euler equations)
	\end{example}
	We simulate the forward step problem \cite{Woodward-Colella1984} of two-dimensional Euler equations \eqref{euler2d} with uniform Mach $3$ flow in a wind tunnel containing a step $[0,3] \times[0,1] \setminus \{[0.6,3] \times [0,0.2]\}$.
	Initially, the wind tunnel is filled with a $\gamma$-law gas, i.e., everywhere satisfies $(\rho,\mu,\nu, P)=(1.4,3,0,1)$.
	Reflective boundary conditions are used along the wall of the tunnel.
	Inflow and outflow boundary conditions are used at the entrance and exit, respectively.
	
\begin{figure}[H]
		\centering
		\centering
		\subfigure[HWENO-R]{
			\includegraphics[width=0.8\textwidth,trim=20 10 20 0,clip]{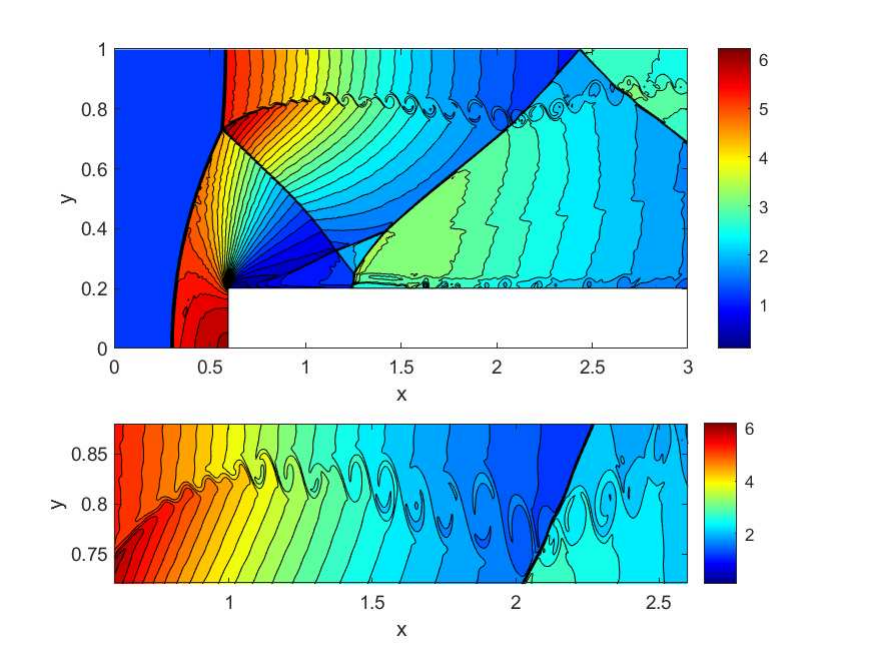}   }
		\subfigure[HWENO-I]{
			\includegraphics[width=0.8\textwidth,trim=20 10 20 0,clip]{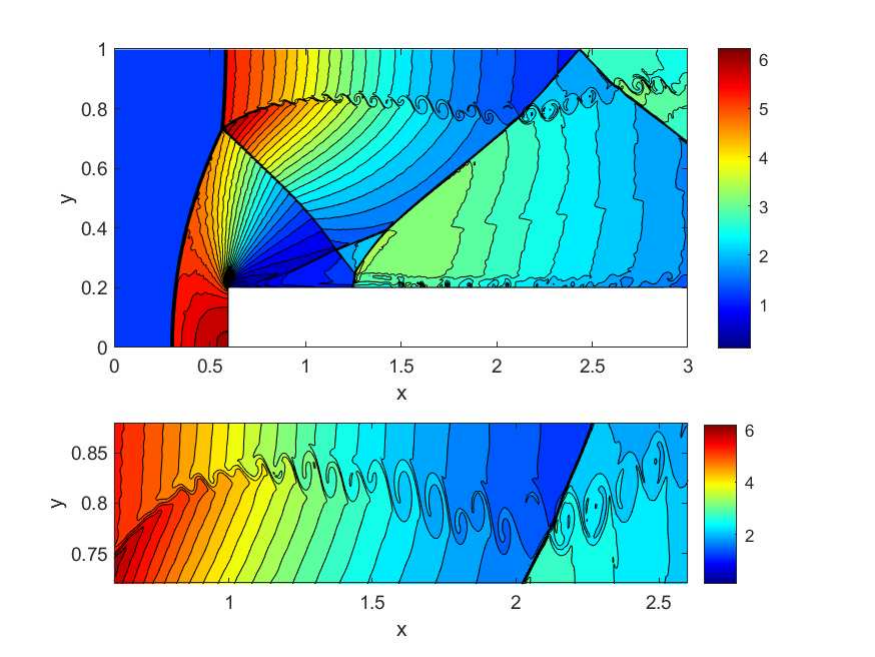}    }
		\caption{{\bf Example} \ref{forwardstep}.
			The density $\rho$ at $t=4$ obtained by the HWENO-R and HWENO-I schemes with $N_x \times N_y=960\times 320$.}\label{Fig:forwardstep}
	\end{figure}

	The density $\rho$ at $t=4$ obtained by the HWENO-R and HWENO-I schemes with $N_x \times N_y=960\times 320$ is shown in Fig. \ref{Fig:forwardstep}.
	The figures show that the results are comparable, and the HWENO-I scheme has a slightly better resolution near the region $[1,2.5]\times[0.75,0.85]$.
	The total CPU time of the HWENO-R scheme and the HWENO-I scheme with $N_x \times N_y=960\times 320$ are $49504.96$s and $35655.17$s, respectively, which shows that the cost by the HWENO-I only about $72\%$ of that by the HWENO-R scheme.
	Again, the HWENO-I scheme has higher efficiency than the HWENO-R scheme.
	
	\begin{example}\label{Mach10shock}
		(Mach 10 shock problem with shock reflection and diffraction of the two-dimensional Euler equations)
	\end{example}
	We consider a Mach $10$ shock problem with shock reflection and diffraction of the two-dimensional Euler equations \eqref{euler2d} over computational
	domain $[0,3]\times[-1,1] \setminus\{[0,1]\times[-1,0]\}$. The initial condition is set as the same as Example \ref{doublemach}, i.e.,
	\begin{equation*}
		(\rho, \mu, \nu, P)=\begin{cases}
			\left(8, 8.25 \cos\left(\frac{\pi}{6}\right),-8.25 \sin\left(\frac{\pi}{6}\right), ~116.5\right),  &\text{if} ~x<\frac{1}{6}+\frac{y}{\sqrt{3}}, \\
			(1.4,~0,~0,~1) & \text { otherwise. }
		\end{cases}
	\end{equation*}
	We use the exact motion of a Mach $10$ shock on the top boundary $[0,3]\times \{1\}$, the exact post-shock condition on $\left[0, \frac{1}{6}\right] \times\{0\}$, and the reflective boundary on both $\left[\frac{1}{6}, 1\right] \times\{0\}$ and $\{1\} \times[-1,0]$. The inflow boundary condition is used on $\{0\}\times [0,1]$ (left) and outflow boundary conditions are used on boundaries $\{3\}\times[-1,1]$ (right) and $[1,3]\times \{-1\}$ (bottom).

\begin{figure}[H]
\centering
\subfigure[HWENO-R]{
\includegraphics[width=0.45\textwidth,trim=10 30 20 50,clip]{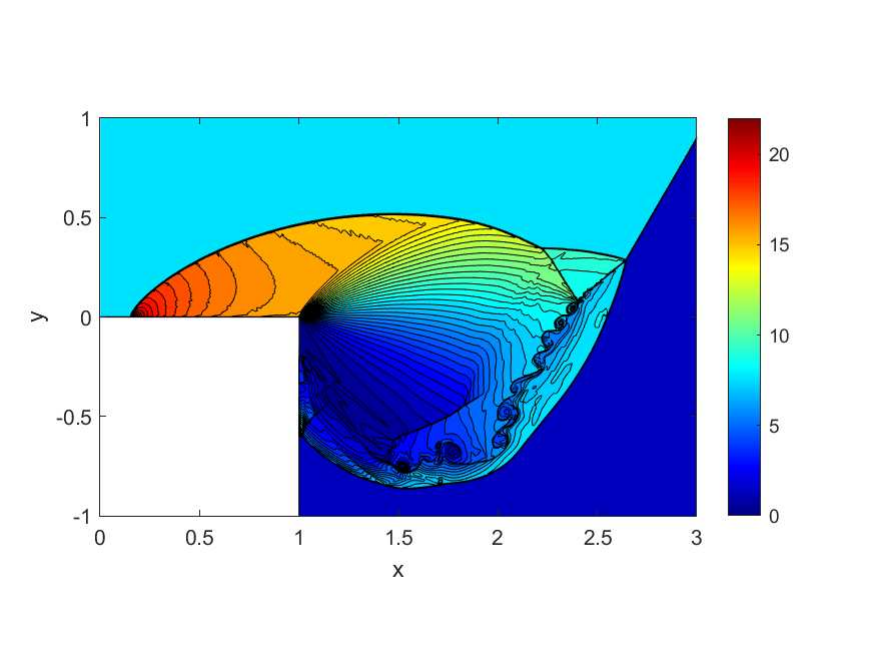}}
\subfigure[HWENO-I]{
\includegraphics[width=0.45\textwidth,trim=10 30 20 50,clip]{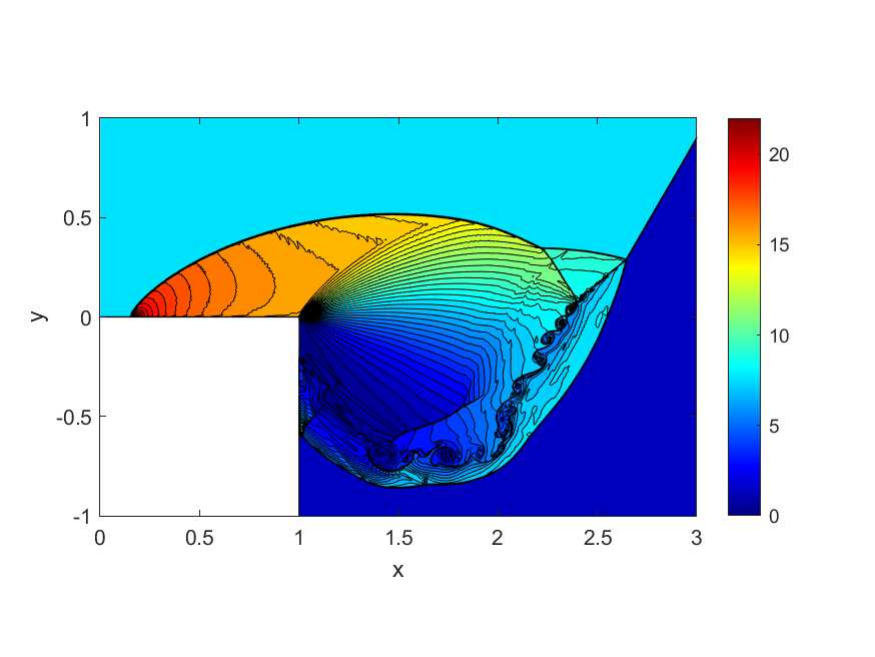}}
\subfigure[HWENO-R at $(1.2, 2.5)\times(0.9,0.2)$]{
\includegraphics[width=0.45\textwidth, trim=10 10 20 10,clip]{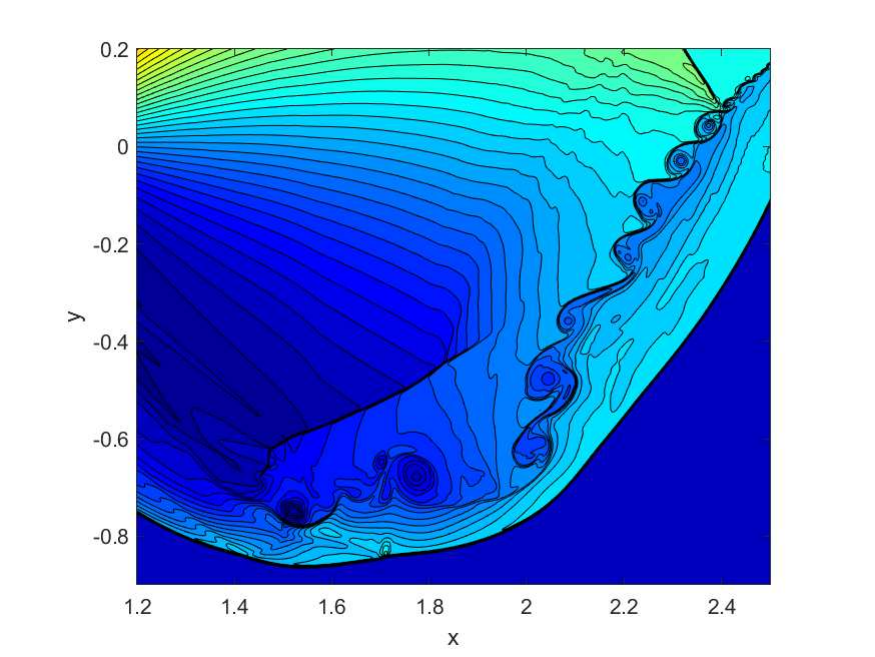}}
\subfigure[HWENO-I at $(1.2, 2.5)\times(0.9,0.2)$]{
\includegraphics[width=0.45\textwidth, trim=10 10 20 10,clip]{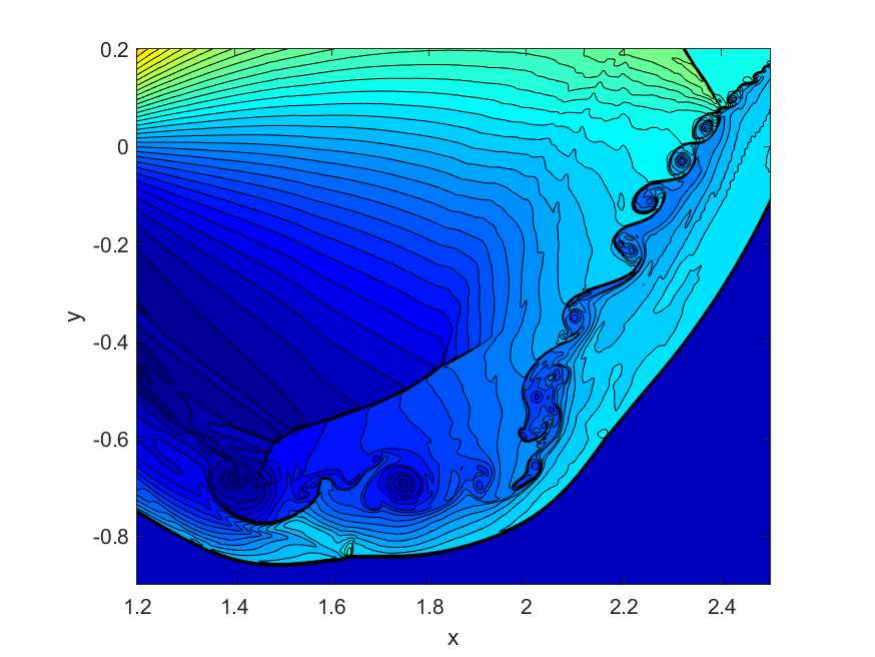}}
\caption{{\bf Example} \ref{Mach10shock}.
{
The density at $t=0.2$ obtained by the the HWENO-R and HWENO-I schemes with $N_x\times N_y=1440 \times 960$.}}
\label{Fig:mach10}
\end{figure}

{
The HWENO-I scheme could work well in simulating the Mach $10$ shock problem with shock reflection and diffraction.
	The density at $t=0.2$ obtained by the HWENO-R and HWENO-I schemes with $N_x\times N_y=1440 \times 960$ (with a parametrized PP flux limiter \cite{XiongQiuXu-PP}) is shown in Fig.~\ref{Fig:mach10}.
	We can see that both schemes have comparable resolution for capturing complex vortex structures.
}

\begin{example}\label{Mach2000}
		(Mach $2000$ astrophysical jet problem of the two-dimensional Euler equations)
	\end{example}
	In this example, we test the Mach $2000$ astrophysical jet problem \cite{Mach80jet} of the two-dimensional Euler equations \eqref{euler2d} with the ratio of specific heats $\gamma = \frac{5}{3}$.
	The computational domain is $[0,1] \times [-0.25,0.25]$.
	The initial condition is set as $(\rho, \mu, \nu, P)=(0.5,0,0,0.4127)$.
	Outflow boundary conditions are posed on the right, top, and bottom, and set
	\begin{equation*}
		(\rho, \mu, \nu, P)=\begin{cases}
			(5,800,0,0.4127), & \text{if}~(x,y)\in\{0\}\times [-0.05,0.05],\\
			(0.5,0,0,0.4127), & \text{if}~(x,y)\in\{0\}\times [-0.25,-0.05] \cup[0.05,0.25],
		\end{cases}
	\end{equation*}
	on the left boundary.
	This is also an extreme problem and many schemes (e.g., \cite{HWENO-R, HWENO-L, WENO-ZhuQiu2016JCP}) fail to simulate this problem without PP limiters.
	
\begin{figure}[H]
		\centering
		\subfigure[density $\rho$ logarithm]{
			\includegraphics[width=0.45\textwidth,trim=10 60 20 50,clip]{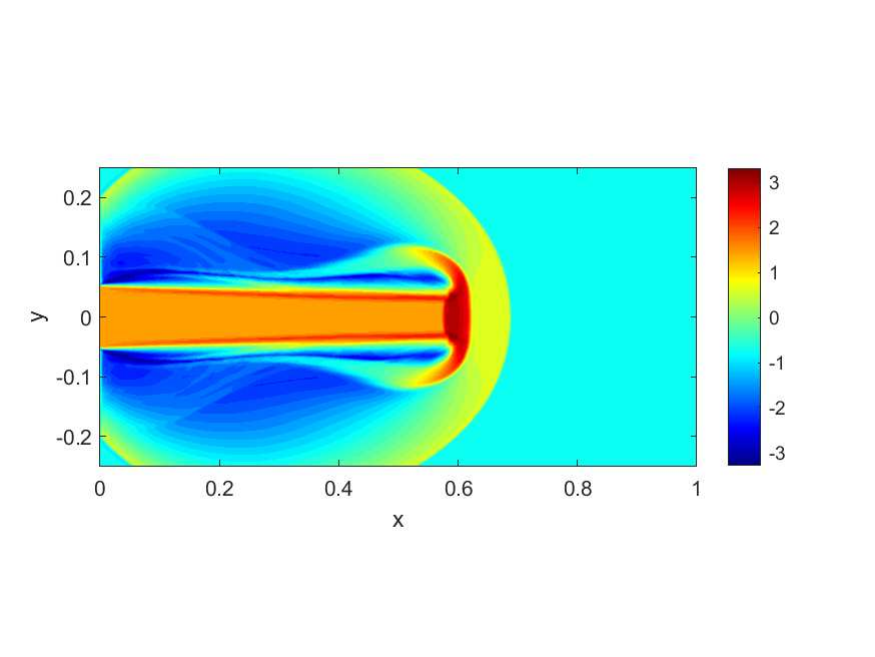}}
		\subfigure[velocity $\mu$ magnitude]{
			\includegraphics[width=0.45\textwidth,trim=10 60 20 50,clip]{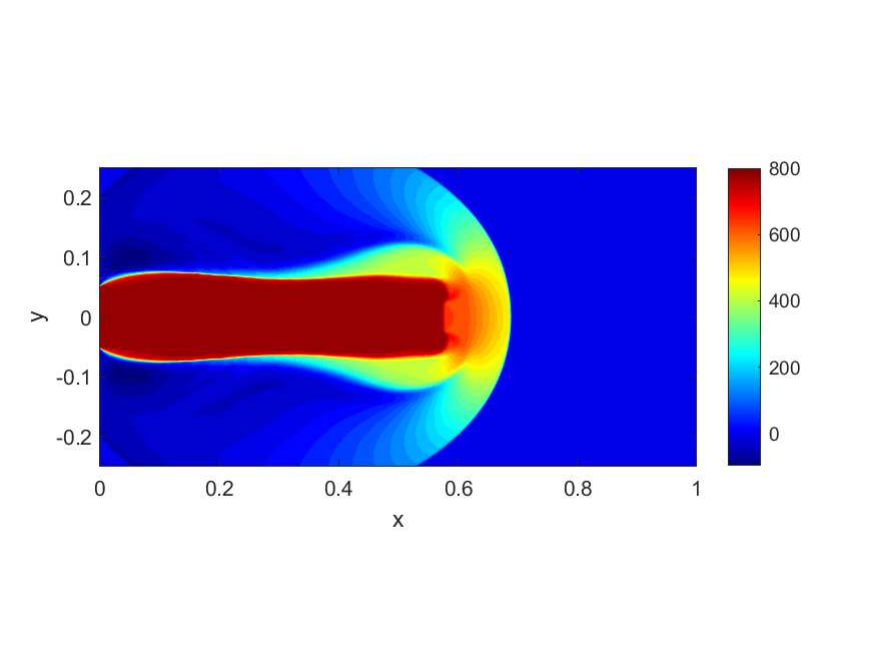}}
		\subfigure[pressure $P$ logarithm]{
			\includegraphics[width=0.45\textwidth,trim=10 60 20 50,clip]{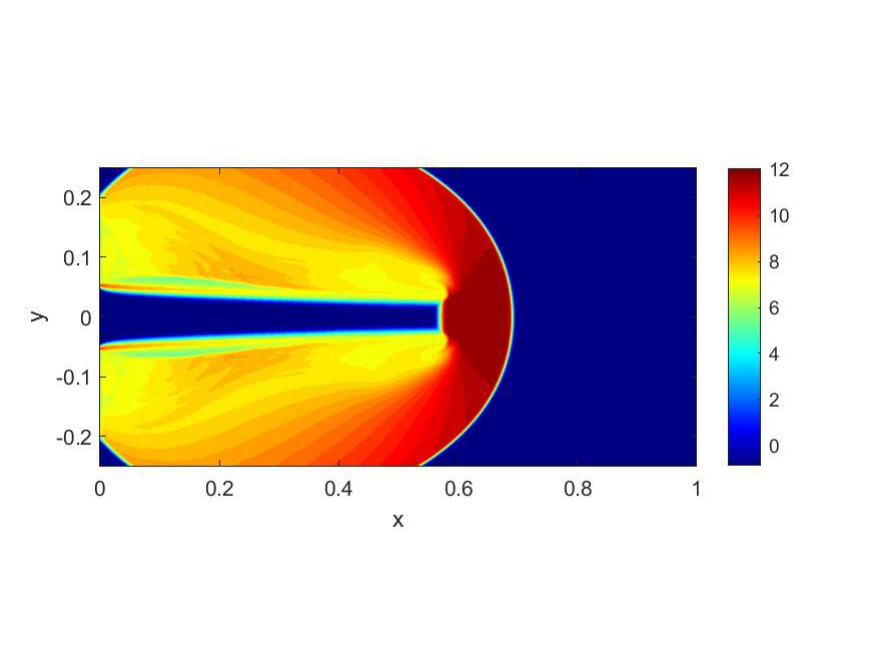}}
		\subfigure[temperature $T$ logarithm]{
			\includegraphics[width=0.45\textwidth,trim=10 60 20 50,clip]{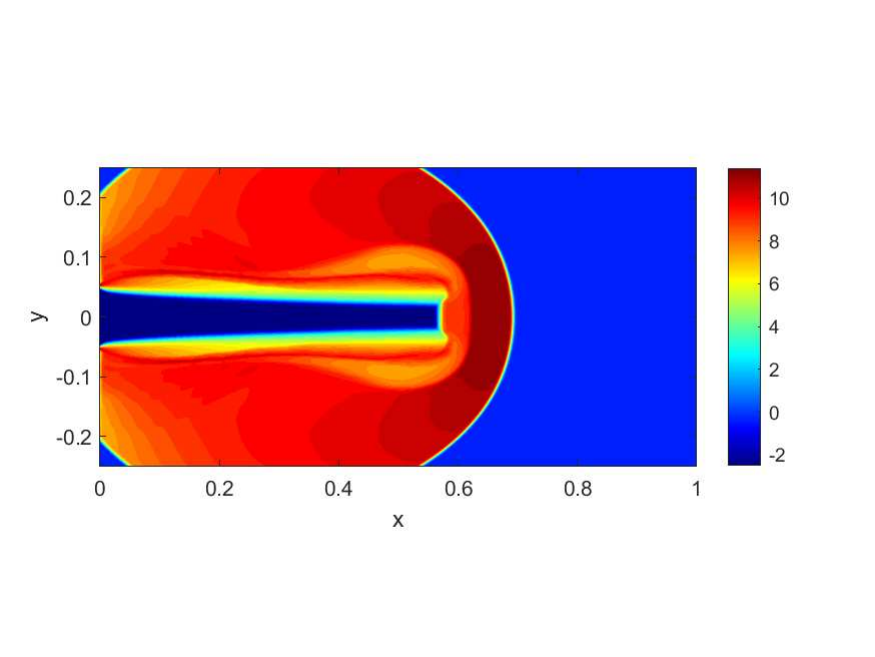}    }
		\caption{{\bf Example} \ref{Mach2000}.
			The density $\rho$, velocity $\mu$, pressure $P$, and temperature $T=P/\rho$ at $t=0.001$ obtained by the HWENO-I scheme with $\gamma_0=0.8$ and $N_x\times N_y=640\times 320$.}
		\label{Fig:mach2000-080}
	\end{figure}

\begin{figure}[H]
		\centering
		\subfigure[density $\rho$ logarithm]{
			\includegraphics[width=0.45\textwidth,trim=10 60 20 50,clip]{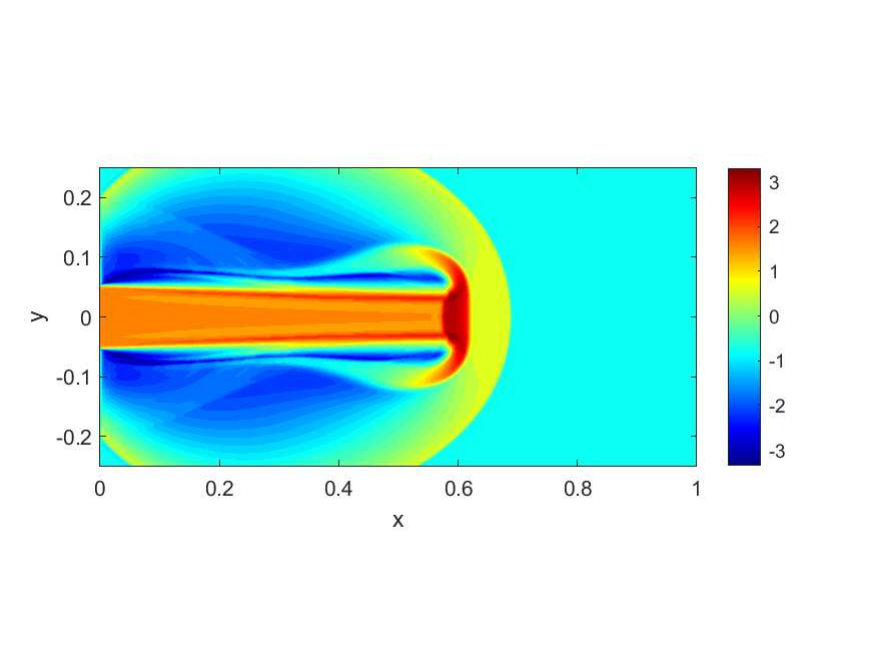}   }
		\subfigure[velocity $\mu$ magnitude]{
			\includegraphics[width=0.45\textwidth,trim=10 60 20 50,clip]{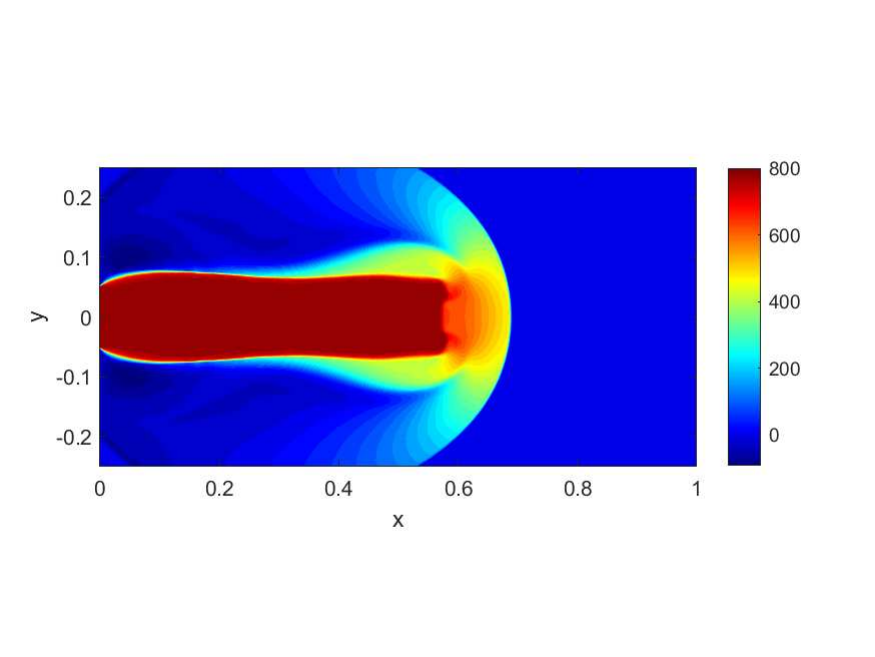}    }
		\subfigure[pressure $P$ logarithm]{
			\includegraphics[width=0.45\textwidth,trim=10 60 20 50,clip]{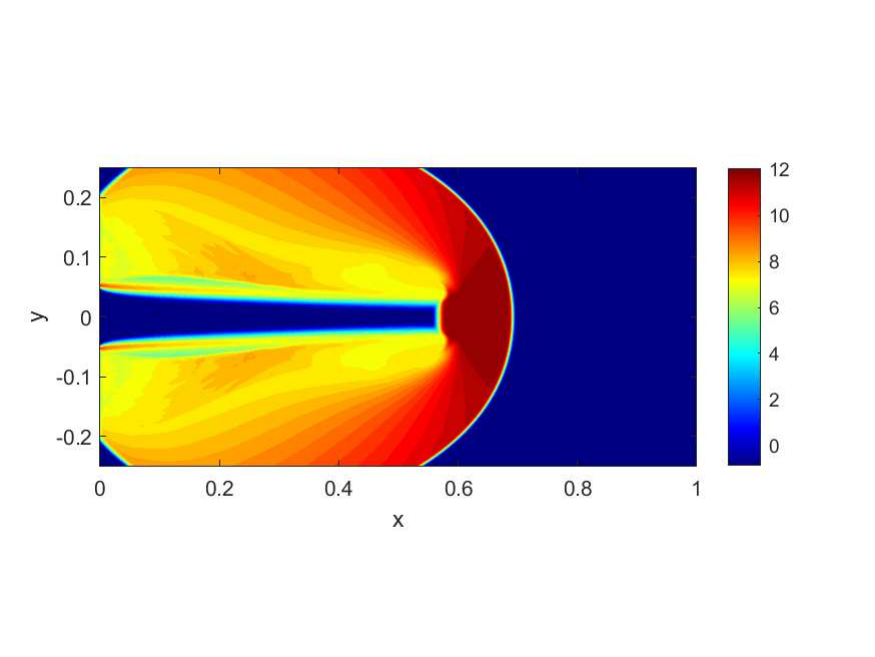}}
		\subfigure[temperature $T$ logarithm]{
			\includegraphics[width=0.45\textwidth,trim=10 60 20 50,clip]{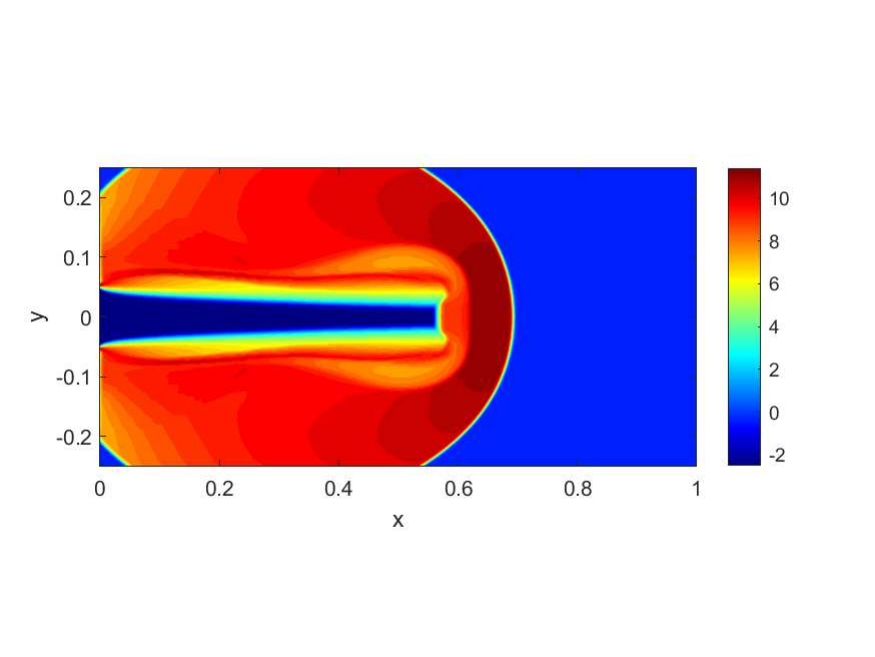}    }
		\caption{{\bf Example} \ref{Mach2000}.
			The density $\rho$, velocity $\mu$, pressure $P$, and temperature $T=P/\rho$ at $t=0.001$ obtained by the HWENO-I scheme with $\gamma_0=1/3$ and $N_x\times N_y=640\times 320$.}
		\label{Fig:mach2000-033}
	\end{figure}

	Interestingly, by trial and error, we find that the proposed HWENO-I scheme can work well with a suitable linear weight $\gamma_0\leq 0.8$ while setting a small initial time step.
	In our computation, the time step $\Delta t$ is initially set as $=10^{-7}$, and then to be chosen according to the CFL condition.
	The density, velocity, pressure, and temperature at $t=0.001$ obtained by the HWENO-I scheme with
	$\gamma_0=0.8$ and $\gamma_0=\frac{1}{3}$ are shown in Fig.~\ref{Fig:mach2000-080} and ~\ref{Fig:mach2000-033}, respectively.
	We can see that the HWENO-I scheme with both $\gamma_0=0.80$ (a relatively big value) and $\gamma_0=\frac{1}{3}$ (a relatively small value) in our implementation work well.

	We guess the reason is that if the time step $\Delta t$ is calculated by CFL condition directly at initial time $t=0$, it may result in an excessively large $\Delta t$ (e.g., about $3.9965\times 10^{-4}$ with $N_x\times N_y=640\times 320$), and this will affect the computational stability.
	
	\begin{example}\label{RTinstability}
		(Rayleigh-Taylor instability problem of the two-dimensional Euler equations with the source terms)
	\end{example}
	In the final example, we simulate the Rayleigh-Taylor instability problem \cite{ShiZhangShu2003} of two-dimensional Euler equations with the source terms
	\begin{equation}\label{euler2d-s}
		\frac{\partial}{\partial t}
		\begin{bmatrix*}[c]
			\rho \\
			\rho \mu \\
			\rho \nu \\
			E
		\end{bmatrix*}
		+
		\frac{\partial}{\partial x}
		\begin{bmatrix*}[c]
			\rho \mu \\
			\rho \mu^2 +P\\
			\rho \mu\nu \\
			\mu(E+P)
		\end{bmatrix*}
		+
		\frac{\partial}{\partial y}
		\begin{bmatrix*}[c]
			\rho \nu \\
			\rho \mu\nu \\
			\rho \nu^2+P\\
			\nu(E+P)
		\end{bmatrix*}
		=\begin{bmatrix*}[c]
			0 \\
			0\\
			g\rho\\
			g\rho \nu
		\end{bmatrix*},
	\end{equation}
	where $g$ is the gravity taken as $1$ in here.
	Rayleigh-Taylor instability happens on an interface between fluids with different densities when acceleration is directed from the heavy fluid to the light fluid.
	The computational domain is $[0,0.25] \times[0,1]$.
	Initially, the interface is at
	$y=\frac{1}{2}$, the heavy fluid with $\rho=2$ is below the interface, the light fluid $\rho=1$ is above the interface with the acceleration in the positive $y$-direction, the pressure $P$ is continuous across the interface, and a small perturbation is given to the $y$-direction fluid speed. Specifically,
	the initial condition is given by
	\begin{equation*}
		(\rho, \mu, \nu, P)= \begin{cases}
			(2,~0,~-0.025 c \cdot \cos (8 \pi x), ~2 y+1), &  \text{if}~(x,y)\in [0,0.25] \times [0,\frac{1}{2}), \\
			(1,~0,~-0.025 c \cdot \cos (8 \pi x),~ y+1.5~), &  \text{if}~(x,y)\in [0,0.25] \times[\frac{1}{2},1],
		\end{cases}
	\end{equation*}
	where $c=\sqrt{\gamma P / \rho}$ is the sound speed with the ratio of specific heats $\gamma = \frac{5}{3}$.
	Reflective boundary conditions are imposed on the left and right, and we set
	\begin{equation*}
		(\rho, \mu, \nu, P)=\begin{cases}
			(2,~0,~0,~1~~), & \text{if}~(x,y)\in[0,0.25]\times \{0\},\\
			(1,~0,~0,~2.5), & \text{if}~(x,y)\in[0,0.25]\times \{1\},
		\end{cases}
	\end{equation*}
	on the top and bottom boundaries.
	
\begin{figure}[H]
		\centering
		\subfigure[HWENO-R]{
			\includegraphics[width=0.9\textwidth,trim=10 10 20 0,clip]{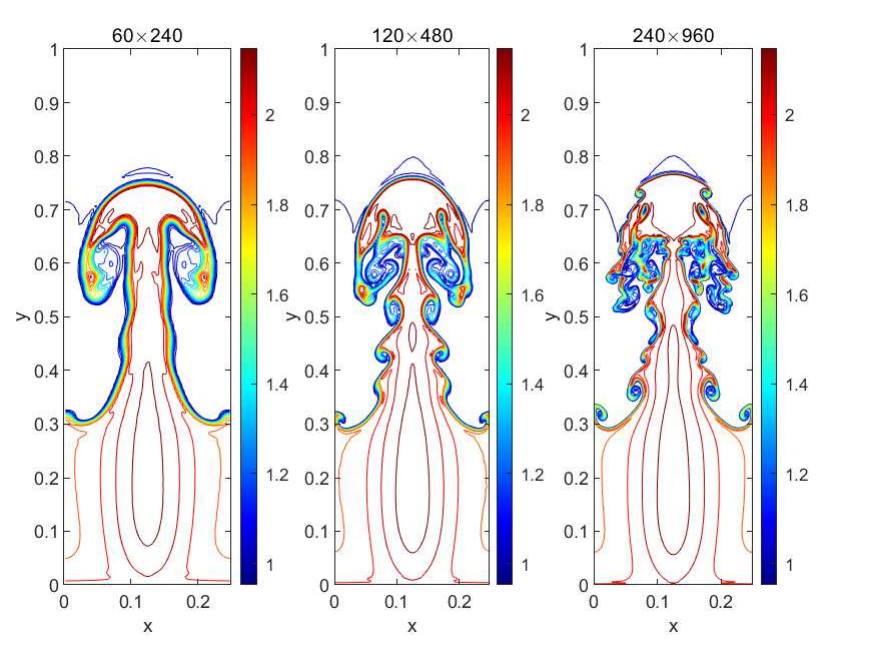}    }
		\subfigure[HWENO-I]{
			\includegraphics[width=0.9\textwidth,trim=10 10 20 0,clip]{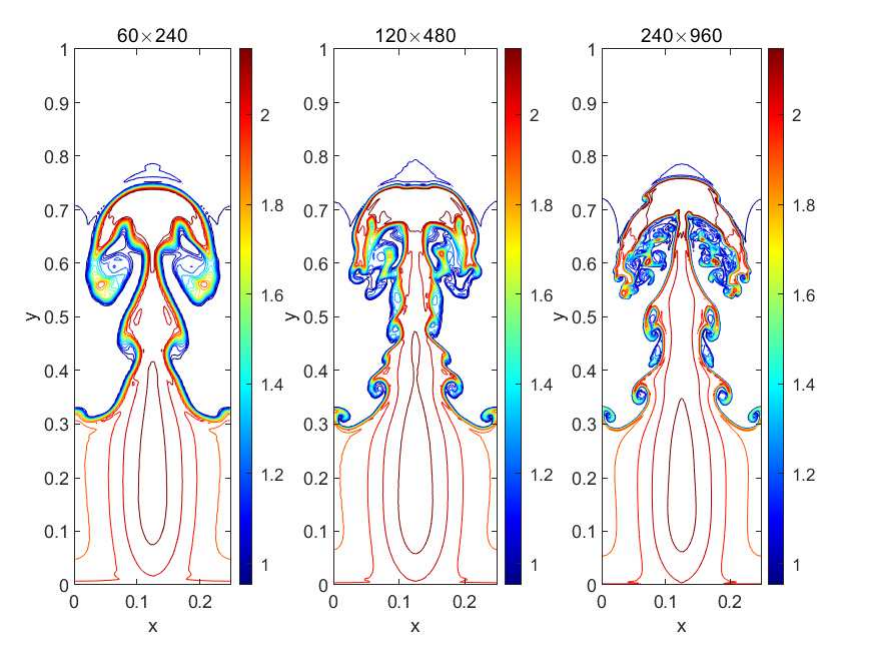}  }
		\caption{{\bf Example} \ref{RTinstability}.
			The density at $t=1.95$ obtained by the HWENO-R and HWENO-I schemes with $N_x\times N_y=60\times 240$, $120\times 480$, and $240\times 960$.}
		\label{Fig:RTinstability-1}
	\end{figure}

	The instability has a fingering nature with bubbles of light fluid rising into the ambient heavy fluid and spikes of heavy fluid falling into the light fluid.
	The density at $t=1.95$ obtained by the HWENO-R and HWENO-I schemes with $N_x\times N_y=60\times 240$, $120\times 480$, and $240\times 960$ are shown in Fig.~\ref{Fig:RTinstability-1}.
	One can see that the HWENO-I scheme can capture complicated solution
	structures well.
	In addition, the CPU times of them are listed in Table.~\ref{tab:RTcpu}, and we can observe that when comparable resolution is obtained by the HWENO-I and HWENO-R schemes, the CPU time needed by the HWENO-I is only about $74-78\%$ of that required by the HWENO-R in our implementation.
	Thus, the HWENO-I scheme is more efficient than the HWENO-R scheme.
	
	\begin{table}[H]
		\caption{{\bf Example} \ref{RTinstability}. CPU time in seconds by the HWENO-R and HWENO-I schemes for the Rayleigh-Taylor instability problem. $t=1.95$.}
		\label{tab:RTcpu}
		\vspace{5pt}
		\centering
		\begin{tabular}{ c c c c  c}
			\toprule
			$N_x\times N_y$ &HWENO-R (s) & HWENO-I (s) & Efficiency (\%) \\ \midrule
			$60\times 240$&302.03&   229.58  & 76.0 \\
			$120\times 480$&2377.70& 1843.07 & 77.5  \\
			$240\times 960$&23563.29& 17574.21& 74.6 \\
			\bottomrule
		\end{tabular}
	\end{table}
	
	\section{Conclusions}
	\label{sec:conclusions}
	In this paper, a simple, efficient, and fifth-order finite difference interpolation-based Hermite WENO (HWENO-I) scheme for the one- and two-dimensional hyperbolic conservation laws has been developed.
	In the reconstruction of numerical fluxes, we do not need to split the flux functions. We first interpolate the values of the solution $u^{\pm}_{i+\frac12}$ and first-order derivative $v^{\pm}_{i+\frac12}$ at $x_{i+\frac12}$ and evaluate the numerical fluxes based on these interpolated values.
	The HWENO interpolation only needs to be applied one time which utilizes the same candidate stencils and Hermite interpolation polynomials in both interpolations of the values $u^{-}_{i+\frac12}$, $v^{-}_{i+\frac12}$, and the modified derivative $\tilde{v}_{i}$ at $x_i$.
	The linear weights in our HWENO-I scheme can be any positive numbers as long as their sum equals to one.
	
	The implementation of the HWENO-I scheme is simpler and more efficient than the previous interpolation-based HWENO scheme \cite{LiuQiu2014JSC-2} which must split the flux for the stability and upwind performance of the approximation for high-order derivative terms.
	The HWENO-I scheme also inherits the advantages of the flux-reconstruction-based HWENO-R scheme \cite{HWENO-R}, including fifth-order accuracy, compact stencils, arbitrary positive linear weights, and high resolution.
	In addition, our HWENO-I scheme needs neither flux splitting nor additional HWENO interpolation for the modified derivative.
	This makes the HWENO-I scheme more efficient than the HWENO-R scheme.

	Various benchmark numerical examples have been simulated to demonstrate the accuracy, efficiency, resolution, and robustness of the proposed HWENO-I scheme.
	We have concluded that when comparable resolution is obtained with the HWENO-I and HWENO-R schemes, the CPU time by the HWENO-I is only about $72-78\%$ of that required by the HWENO-R scheme.
	The choice of linear weights generally does not affect the accuracy of the HWENO-R and HWENO-I schemes, however, it could have an impact on the resolution of the scheme.
	The smaller $\gamma_0$ used in HWENO interpolation/reconstruction could bring better robustness. Meanwhile, the  resolution will also deteriorate.
	It is interesting that, as $\gamma_0$ decreases, the increase in numerical dissipation of the HWENO-I scheme seems slower than that of the HWENO-R scheme, that is to say, the HWENO-I scheme has better resolution than the HWENO-R scheme (cf. Example \ref{blast-1d}).
	
	\section*{Acknowledgements}
	Min Zhang is partially supported by the National Natural Science Foundation of China (Grant Number: 12301493).
	The computational resources were partially supported by High-performance Computing Platform of Peking University.


\begin{thebibliography}{10}
		
		\bibitem{Capdeville}
		G. Capdeville,
		A Hermite upwind WENO scheme for solving hyperbolic conservation laws, {\em J. Comput. Phys.},
		227 (2008), 2430-2454.
		
		\bibitem{HWENO-U-arxiv}
		C. Fan, J. Qiu and Z. Zhao,
		A moment-based Hermite WENO scheme with unified stencils for hyperbolic conservation laws,
		{\em arxiv:2402.03074.}
		
		\bibitem{HWENO-R}
		C. Fan, Z. Zhao, T. Xiong, and J. Qiu,
		A robust fifth order finite difference Hermite WENO scheme for compressible Euler equations,
		{\em Comput. Methods Appl. Mech. Engrg.}, 412 (2023), 116077.
		
		\bibitem{AWENO-9}
		Z. Gao, L.L. Fang, B.S. Wang, Y. Wang, and W.S. Don,
		Seventh and ninth orders characteristic-wise alternative WENO finite difference schemes for hyperbolic conservation laws,
		{\em Comput.} \& {\em Fluids}, 202 (2020), 104519.
		
		\bibitem{Mach80jet}
		Y. Ha, C.L. Gardner, A. Gelb, and C.-W. Shu,
		Numerical simulation of high Mach number astrophysical jets with radiative cooling,
		{\em J. Sci. Comput.}, 24 (2005), {29-44}.
		
		\bibitem{ENO-1}
		A. Harten,  B. Engquist, S. Osher, and S. R. Chakravarthy,
		Uniformly high order accurate essentially non-oscillatory schemes, III,
		{\em J. Comput. Phys.}, 71 (1987), 231-303.
		
		\bibitem{ENO-2}
		A. Harten and S. Osher,
		Uniformly high-order accurate nonoscillatory schemes I,
		{\em SIAM J. Numer. Anal.}, 24 (1987), 279-309.
		
		\bibitem{Hu-PP2013}
		X.Y. Hu, N.A. Adams, and C.-W. Shu
		Positivity-preserving method for high-order conservative schemes
		solving compressible Euler equations,
		{\em J. Comput. Phys.}, 242 (2013), {169-180}.
		
		\bibitem{WENOJS1996}
		G.-S. Jiang and C.-W. Shu,
		Efficient implementation of weighted ENO schemes,
		{\em J. Comput. Phys.}, 126 (1996), 202-228.
		
		\bibitem{JiangShuZhang2013}
		Y. Jiang, C.-W. Shu, and M. Zhang,
		An alternative formulation of finite difference weighted ENO schemes with Lax--Wendroff time discretization for conservation laws,
		{\em SIAM J. Sci. Comput.}, 35 (2013), A1137-A1160.
		
		
		\bibitem{LiMRHW1}
		J. Li, C.-W. Shu and J. Qiu,
		Multi-resolution HWENO schemes for hyperbolic conservation laws,
		{\em J. Comput. Phys.}, 446 (2021), 110653.
		
		\bibitem{Linde-Roe-1997}
		T. Linde and  P.L. Roe,
		Robust Euler codes. In: 13th Computational Fluid Dynamics Conference, AIAA-97-2098, AIAA Inc., Reno, Nevada (1997).
		
		
		
		\bibitem{LiuQiu2014JSC}
		H. Liu and J. Qiu,
		Finite difference Hermite WENO Schemes
		for hyperbolic conservation laws,
		{\em J. Sci. Comput.},  63 (2015), {548-572}.
		
		\bibitem{LiuQiu2014JSC-2}
		H. Liu and J. Qiu,
		Finite difference Hermite WENO Schemes for conservation Laws, II: An alternative approach,
		{\em J. Sci. Comput.},  66 (2016), {598-624}.
		
		\bibitem{WENO1994}
		X.D. Liu, S. Osher, and T. Chan,
		Weighted essentially non-oscillatory schemes,
		{\em J. Comput. Phys.}, 115 (1994), 200-212.
		
		\bibitem{HWENO-5}
		Z. Ma and S.-P. Wu, HWENO schemes based on compact difference for hyperbolic conservation laws,
		{\em J. Sci. Comput.}, 76 (2018), 1301-1325.
		
		\bibitem{HWENO-1}
		J. Qiu and C.-W. Shu,
		Hermite WENO schemes and their application as limiters for Runge-Kutta discontinuous Galerkin method: one-dimensional case,
		{\em J. Comput. Phys.}, 193 (2004), 115-135.
		
		\bibitem{HWENO-2}
		J. Qiu and C.-W. Shu,
		Hermite WENO schemes and their application as limiters for Runge-Kutta discontinuous Galerkin method II: Two-dimensional case,
		{\em Comput.} \& {\em Fluids}, 34 (2005), 642-663.
		
		\bibitem{Sedov1959}
		L.I. Sedov,
		Similarity and Dimensional Methods in Mechanics,
		Academic Press, New York (1959).
		
		\bibitem{ShiZhangShu2003}
		J. Shi, Y.-T. Zhang, and C.-W. Shu,
		Resolution of high order WENO schemes for complicated flow structures
		{\em J. Comput. Phys.}, 186 (2003), {690-696}.
		
		\bibitem{WENO-Shu2009}
		C.-W. Shu,
		High order weighted essentially nonoscillatory schemes for convection dominated problems,
		{\em SIAM Review}, 51 (2009), 82-126.
		
		\bibitem{WENO-Shu2020}
		C.-W. Shu,
		Essentially non-oscillatory and weighted essentially non-oscillatory schemes,
		{\em Acta Numerica}, 29 (2020), 701-762.
		
		\bibitem{Shu-Osher1988}
		C.-W. Shu and S. Osher,
		Efficient implementation of essentially non-oscillatory shock-capturing schemes,
		{\em J. Comput. Phys.}, 77 (1988), 439-471.
		
		\bibitem{Shu-Osher1989}
		C.-W. Shu and S. Osher,
		Efficient implementation of essentially non-oscillatory shock-capturing schemes, II,
		{\em J. Comput. Phys.}, 83 (1989), {32-78}.
		
		\bibitem{Tao2015JCP}
		Z. Tao, F. Li, and J. Qiu,
		High-order central Hermite WENO schemes on staggered meshes for hyperbolic conservation laws,
		{\em J. Comput. Phys.}, 281 (2015), {148-176}.
		
		\bibitem{Tao2016JCP}
		Z. Tao, F. Li, and J. Qiu,
		High-order central Hermite WENO schemes: dimension-by-dimension moment-based reconstructions,
		{\em J. Comput. Phys.}, 318 (2016), 222-251.
		
		\bibitem{Tao2024JCP}
		Z. Tao, J. Zhang, J. Zhu, and J. Qiu, High-order multi-resolution central Hermite WENO schemes for hyperbolic conservation laws,
		{\em J. Sci. Comput.}, 99 (2024): 40.
		
		\bibitem{AWENO-1}
		B.S. Wang, P. Li, Z. Gao, and W.S. Don,
		An improved fifth order alternative WENO-Z finite difference scheme for hyperbolic conservation laws,
		{\em J. Comput. Phys.}, 374 (2018), 469-477.
		
		\bibitem{AWENO-MR2021}
		Z. Wang, J. Zhu, Y. Yang, and N. Zhao,
		A new fifth-order alternative finite difference multi-resolution WENO
		scheme for solving compressible flow,
		{\em Comput. Methods Appl. Mech. Engrg.},  382 (2021), 113853.
		
		\bibitem{Woodward-Colella1984}
		P. Woodward and P. Colella,
		The numerical simulation of two-dimensional fluid flow with strong shocks,
		{\em J. Comput. Phys.}, 54 (1984), 115-173.
		
		\bibitem{XiongQiuXu-PP}
		T. Xiong, J.-M. Qiu, and Z. Xu,
		Parametrized positivity preserving flux limiters for the high order finite difference WENO scheme solving compressible Euler equations,
		{\em J. Sci. Comput.}, 67 (2016), 1066-1088.
		
		\bibitem{7thHWENO}
		Y.H. Zahran and  A.H. Abdalla, Seventh order Hermite WENO scheme for hyperbolic conservation laws,
		{\em Comput.} \& {\em Fluids}, (2016) 66-80.
		
		\bibitem{HWENO-L}
		M. Zhang and Z. Zhao,
		A fifth-order finite difference HWENO scheme combined with limiter for hyperbolic conservation laws,
		{\em J. Comput. Phys.}, 472 (2023), 111676.
		
		\bibitem{HWENO-a}
		Z. Zhao and J. Qiu,
		A Hermite WENO scheme with artificial linear weights for hyperbolic conservation laws,
		{\em J. Comput. Phys.}, 417 (2020), 109583.
		
		\bibitem{HWENO-M}
		Z. Zhao, Y.-T. Zhang, and J. Qiu,
		A modified fifth order finite difference Hermite WENO scheme for hyperbolic conservation laws,
		{\em J. Sci. Comput.}, 85 (2020), 29.
		
		\bibitem{HWENO-3}
		J. Zhu and J. Qiu,
		A class of the fourth order finite volume Hermite weighted essentially non-oscillatory schemes,
		{\em Sci. China Ser. A Math.}, 51 (2008), 1549-1560.
		
		\bibitem{WENO-ZhuQiu2016JCP}
		J. Zhu and J. Qiu,
		A new fifth order finite difference WENO scheme for solving hyperbolic conservation laws,
		{\em J. Comput. Phys.}, 318 (2016), 110-121.
		
		\bibitem{WENO-ZhuShu2018JCP}
		J. Zhu and C.-W. Shu, A new type of multi-resolution WENO schemes with increasingly higher order of accuracy,
		{\em J. Comput. Phys.},  375 (2018), 659-683.
		
		
		
	\end{thebibliography}
\end{document}